\documentclass[11pt]{article}
\usepackage[utf8]{inputenc}
\usepackage{amsmath,amssymb,epsfig,bbm}
\usepackage{stmaryrd,mathabx}
\usepackage{comment}
\usepackage{color}
\usepackage[T1]{fontenc}

\usepackage{interval}

\usepackage[textsize=tiny,textwidth=20mm]{todonotes}
\usepackage{enumitem}
\usepackage{varwidth}
\setlist{nolistsep}
\usepackage{hyperref}

\usepackage{amsthm}

\newtheorem{defi}{Definition}
\newtheorem{prop}[defi]{Proposition}
\newtheorem{theo}[defi]{Theorem}
\newtheorem{theofr}[defi]{Théorème}
\newtheorem{conj}[defi]{Conjecture}
\newtheorem{lemm}[defi]{Lemma}
\newtheorem{lemmfr}[defi]{Lemme}
\newtheorem{coro}[defi]{Corollary}

\theoremstyle{definition}
\newtheorem{rema}[defi]{Remark}
\newtheorem{exem}[defi]{Example}
\newtheorem{exems}[defi]{Examples}

\newcommand{\bdefi}{\begin{defi}}
\newcommand{\edefi}{\end{defi}}
\newcommand{\bprop}{\begin{prop}}
\newcommand{\eprop}{\end{prop}}
\newcommand{\btheo}{\begin{theo}}
\newcommand{\etheo}{\end{theo}}
\newcommand{\btheofr}{\begin{theofr}}
\newcommand{\etheofr}{\end{theofr}}
\newcommand{\blemm}{\begin{lemm}}
\newcommand{\elemm}{\end{lemm}}
\newcommand{\blemmfr}{\begin{lemmfr}}
\newcommand{\elemmfr}{\end{lemmfr}}
\newcommand{\brema}{\begin{rema}}
\newcommand{\erema}{\end{rema}}
\newcommand{\bexer}{\begin{exem}}
\newcommand{\eexer}{\end{exem}}
\newcommand{\bexems}{\begin{exems}}
\newcommand{\eexems}{\end{exems}}
\newcommand{\bconj}{\begin{conj}}
\newcommand{\econj}{\end{conj}}
\newcommand{\bcoro}{\begin{coro}}
\newcommand{\ecoro}{\end{coro}}
\newcommand{\dem}{\noindent{\bf Proof. }}

\usepackage{mathrsfs}
\renewcommand\mathcal{\mathscr}

\newcommand{\I}{{\cal I}}

\newcommand{\N}{{\cal N}}
\newcommand{\OOO}{{\cal O}}

\newcommand{\maths}[1]{{\mathbb #1}}  

\newcommand{\CC}{\maths{C}}

\newcommand{\HH}{\maths{H}}

\newcommand{\KK}{\maths{K}}

\newcommand{\NN}{\maths{N}}

\newcommand{\PP}{\maths{P}}
\newcommand{\QQ}{\maths{Q}}
\newcommand{\RR}{\maths{R}}
\newcommand{\SSS}{\maths{S}}

\newcommand{\ZZ}{\maths{Z}}

\newcommand{\aaa}{{\mathfrak a}}

\newcommand{\ra}{\rightarrow}
\newcommand{\bs}{\backslash}

\newcommand{\ov}[1]{{\overline #1}} 
\newcommand{\wt}[1]{{\widetilde{#1}}}

\newcommand{\ga}{\gamma}
\newcommand{\Ga}{\Gamma}

 \usepackage{mathtools}

\newcommand{\fn}{\mathfrak N}

\newcommand{\cqfd}{\hfill$\Box$}

\newcommand{\bigO}{\operatorname{O}}

\newcommand{\card}{\operatorname{Card}}

\newcommand{\Det}{\operatorname{Det}}

\newcommand{\disc}{D}

\newcommand{\genfar}{$K$-Fa\-rey }

\newcommand{\icg}[1]{\I_{#1}}
\newcommand{\id}{\operatorname{id}}
\renewcommand{\Im}{\operatorname{Im}}

\newcommand{\Nr}{\operatorname{{\tt N}}}

\newcommand{\radix}{f_K}
\renewcommand{\Re}{\operatorname{Re}}

\newcommand{\ssm}{\!\smallsetminus\!}

\newcommand{\sym}{\mathfrak{s}}
\newcommand{\Sym}{\mathfrak{S}_{\Ga_\QQ}}
\newcommand{\Symrec}{\mathfrak{S}^{\rm rec}_{\Ga_\QQ}}

\newcommand{\Vol}{\operatorname{Vol}}
\newcommand{\vol}{\operatorname{vol}}

\newcommand{\hdr}{{\HH}^2_\RR}
\newcommand{\htr}{{\HH}^3_\RR}
\newcommand{\hcr}{{\HH}^5_\RR}
\newcommand{\hnr}{{\HH}^n_\RR}

\newcommand{\PSL}{\operatorname{PSL}}
\newcommand{\SL}{\operatorname{SL}}
\newcommand{\GL}{\operatorname{GL}}

\newcommand{\PSLOK}{\operatorname{PSL}_{2}(\OOO_K)}

\newcommand{\PSLR}{\operatorname{PSL}_{2}(\RR)}

\newcommand{\PSLZ}{\operatorname{PSL}_{2}(\ZZ)}

\newcommand{\tr}{\operatorname{\tt tr}}
\newcommand{\n}{\operatorname{\tt n}}

\newcounter{fig}

\usepackage[percent]{overpic}
\usepackage{pict2e}

\title{Farey neighbours, modular symbols and divergent geodesics}
\author{Jouni Parkkonen \and Fr\'ed\'eric Paulin}
\date{\today}

\begin{document}
\bibliographystyle{../alphanum}
\maketitle

\begin{abstract} 
We give effective asymptotic counting results for pairs of Farey
neighbours and for modular symbols in $\QQ$, in imaginary quadratic
number fields and in definite quaternion algebras over $\QQ$, using
the distribution of common perpendiculars between Margulis cusp
neighbourhoods and divergent geodesics in hyperbolic manifolds. We
describe the tangency properties of the canonical Margulis cusp
neighbourhoods in Bianchi hyperbolic $3$-orbifolds.
\end{abstract}

\section{Introduction}
In this paper, we present effective asymptotic counting results for
pairs of Farey neighbours in $\QQ$, in imaginary quadratic number
fields and in definite quaternion algebras over $\QQ$, when the lower
bound on the distances between the Farey neighbours shrinks to
$0$. These results appear to be new even in the classical rational
case. They are contributions to the study of the distribution of pairs
of the well known Farey fractions and their generalisations, see for
instance \cite{Hall70, HalTen84, Haynes04, BocZah05, Marklof10Inv,
  Marklof10AMS, Marklof13, Athreya16, Heersink16, BocSis22, Lutsko22,
  ParPau24, Sayous25}, and of modular symbols of Shimura, Eichler,
Birch, Manin, see  for instance \cite{Manin72, Cremona84, Manin09,
  McMullen21}.

Let $K$ be either $\QQ$ or an imaginary quadratic number field, with
ring of integers $\OOO_K$, discriminant $D_K$, class number $h_K$ and
Dedekind zeta function $\zeta_K$. Recall that two elements $\alpha,
\beta\in \PP^1(K)=K\cup\{\infty\}$ are {\em Farey neighbours}
if there exists $p,q,r,s\in\OOO_K$ with $\alpha=\frac pq$, $\beta=
\frac rs$ and
\begin{equation}\label{eq:defiFareyNeigh}
  |\;ps-qr\,|=1
\end{equation}
or, equivalently, $ps-qr\in\OOO_K^\times$. The Diophantine equation
$ps-qr=1$ with integral unknowns $p,q,r,s$ is called the {\it gcd
  equation}. See for instance \cite{HorNev23} for other distribution
results of solutions to the gcd equation, and for instance
\cite{Schmidt68, Duke03, EllMicVen13, AkaEinSha16Inv, HorKar23} for
higher dimensional generalisations pioneered by Linnik and Maass.

When the class number $h_K$ of $K$ is greater than $1$, there are
infinitely many elements of $\PP^1(K)$ that do not have Farey
neighbours. In Section \ref{sec:cuspverymax}, we discuss a notion of
generalized Farey neighbours due to Bestvina-Savin \cite{BesSav12},
that geometrically corresponds to the tangency in the real hyperbolic
$3$-space $\htr$ of Mendoza's canonical horoballs \cite{Mendoza80}
centered at two points of $\PP^1(K)$.  We prove in Theorem
\ref{theo:bianverymax} the existence of generalized Farey neighbours
of every element of $\PP^1(K)$, extending \cite[Theo.~2]{Ford38} when
$K=\QQ$.

Let $\fn_K$ be the set of unordered pairs of Farey neighbours in $K$.
The additive group $\OOO_K$ acts by simultaneous translations on the
set $\fn_K$.

\btheo\label{theo:realfareycount}
\hypertarget{mainintro1}{(1)} As $\epsilon>0$ tends to $0$, we
have
\[
\card\big\{ \{\alpha,\beta\}\in \fn_{\QQ} :
\alpha,\beta\in\interval 01\,,\ |\beta-\alpha| \ge \epsilon\big\}
=-\frac 6{\pi^2}\,\frac{\ln \epsilon}\epsilon+\bigO(\epsilon^{-1})\,.
\]
\hypertarget{mainintro2}{(2)} If $K$ is an imaginary quadratic number
field, then as $\epsilon>0$ tends to $0$, we have
\[
\card\big(\OOO_{K}\bs\big\{ \{\alpha,\beta\}\in \fn_K :
|\beta-\alpha| \ge \epsilon\big\}\big) =
\frac {4\,\pi}{|\OOO_K^\times|\;D_K\;\zeta_K(2)}\,
\frac{\ln\epsilon}{\epsilon^2}+\bigO(\epsilon^{-2})\, .
\]
\etheo

We refer to Section \ref{sec:quaternion} (and in particular to Theorem
\ref{theo:quaternionfareycount}) for analogous results on the
effective asymptotic counting function of pairs of Farey neighbours in
definite rational quaternion algebras, noting that their definition
(see Equation \eqref{eq:defiQuatFarNei}) needs to address the
noncommutativity issue of the quaternion algebra.  We have the following
reformulation of the first claim of Theorem \ref{theo:realfareycount} in terms of the solutions of the gcd equation.

\bcoro\label{cor:farey} As $N\in\NN\ssm\{0\}$ tends to $+\infty$, we have
\[
  \card\Big\{ (p,q,r,s)\in\ZZ^4:
  \begin{array}{cc} ps-qr=1\,,\ 0<qs\le N\\
  0\le p\le q\,,\ 0\le r\le s\end{array} \Big\}
  =\frac {12}{\pi^2}\,N\log N+\bigO(N)\,.
\]
\ecoro

The two above results follow (see Subsection \ref{subsec:proofintro})
from the more general Theorem \ref{theo:rationalcount} when $K=\QQ$
and Theorem \ref{theo:iqcount} otherwise, that include versions with
congruences and cover the more general case of modular symbols.
Similarly, Theorem \ref{theo:quaternionfareycount} in the quaternionic
case follows from the more general Theorem \ref{theo:quatercount}.

For every $n\in\NN$, let $\hnr$ be the upper halfspace model of the
real hyperbolic space.  A {\em Farey arc} in $\hdr$ is a hyperbolic
geodesic line in $\hdr$ whose pair of points at infinity is a pair of
rational Farey neighbours. We prove Theorem \ref{theo:rationalcount}
by relating its counting function to the counting function of common
perpendicular geodesic arcs from the horoball $B_\infty=\{z\in\hdr:\Im
z\ge 1\}$ to the Farey arcs.  This correspondence allows us to use the
recent geometric counting results of \cite{ParPau25Ingham} in the
proof.  The arguments for Theorems \ref{theo:iqcount} and
\ref{theo:quaternionfareycount} are similar, using geodesic lines in
$\htr$ and $\HH^5_\RR$ respectively.
 
\begin{center}
\begin{overpic}[width=0.95\textwidth,tics=10]{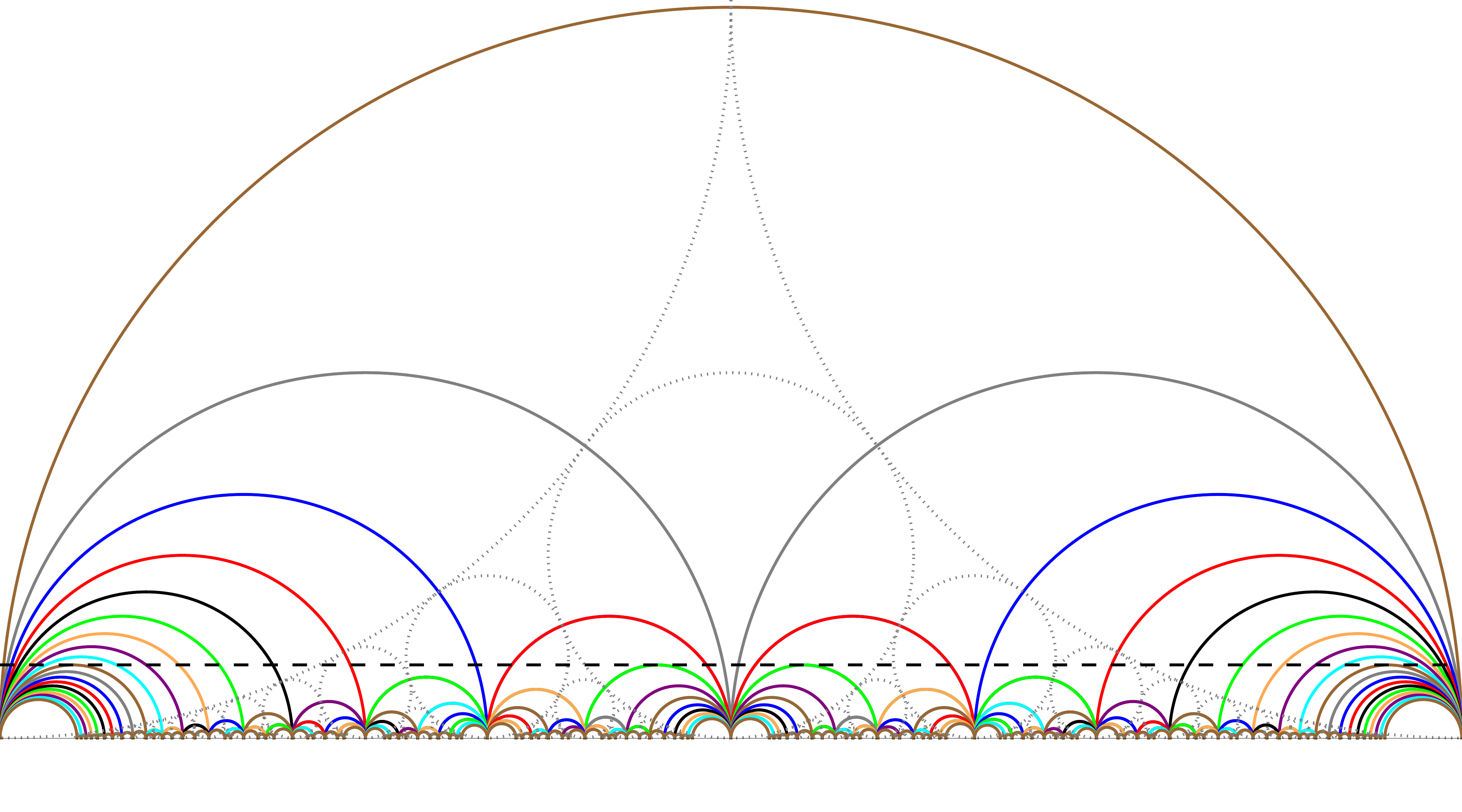} 
  \put (-1,2) {$\frac 01$}
  \put (3.8,2) {$\frac 1{19}$}
  \put (19.1,2) {$\frac 15$}
  \put (24.1,2) {$\frac 14$}
  \put (32.5,2) {$\frac 13$}
  \put (39.1,2) {$\frac 25$}
  \put (49.1,2) {$\frac 12$}
  \put (59.2,2) {$\frac 35$}
  \put (65.8,2) {$\frac 23$}
  \put (74.1,2) {$\frac 34$}
  \put (79.2,2) {$\frac 45$}
  \put (94.2,2) {$\frac {18}{19}$}
  \put (99,2) {$\frac 11$}
\end{overpic}
\end{center}
 
The figure above shows the Farey arcs with endpoints $\frac pq$ and
$\frac rs$ in $\interval 01$ that have denominators $1\le q,s\le 19$
in reduced forms.  The dotted circles show the $\PSLZ$-orbit of the
horoball $B_\infty$.  The horizontal dashed line shows the points in
$\hdr$ at hyperbolic distance $\ln 20$ from $B_\infty$. It meets $23$
Farey arcs with endpoints in $\interval 01$.
\medskip

\noindent{\small {\it Acknowledgements :} We thank Kevin Destagnol
  for his help with Theorem \ref{theo:bianverymax}.}

\section{Counting common perpendiculars to divergent geodesics}
\label{sec:commperp}

For every $n\in\NN\ssm\{0,1\}$, let
\[
\hnr=\Big(\,\{x=(x_1,\dots,x_n)\in\RR^n:\,x_n>0\},
\; ds^2= \frac{dx_1^2+\dots+dx_n^2}{x_n^2}\,\Big)
\]
be the upper halfspace model of the $n$-dimensional real hyperbolic
space with constant sectional curvature $-1$.  Let $\Ga$ be a discrete
group of isometries of $\hnr$ such that $M=\Ga\bs\hnr$ is a finite
volume complete noncompact real hyperbolic orbifold.  Assume that
$\infty$ is a parabolic fixed point of $\Ga$, so that the image of the
horoball $B_\infty=\{x\in\hnr:\,x_n\geq 1\}$ in $M$ is a properly
immersed closed locally convex subset of $M$.

A locally geodesic line $\ell: \RR\to M$ that is a proper mapping is a
{\em divergent geodesic} in $M$.  We denote by $m(\ell(\RR))$ the
cardinality of the orbifold pointwise stabiliser in $\Ga$ of the image
$\ell(\RR)$ of $\ell$.  A locally geodesic line $\ell$ in $M$ (or its
image) is {\em weakly reciprocal} if it has a lift $\wt\ell:\RR\ra
\hnr$ such that an element of $\Ga$ interchanges the two endpoints at
infinity of the geodesic line $\wt\ell$. We say\footnote{As in
\cite{ParPau25Ingham}, see also \cite{Sarnak07}.}  that $\ell$ (or its
image) is {\em reciprocal} if there is such an element of order $2$.
Let $\iota_{\rm rec}(\ell(\RR))=1$ if $\ell$ is weakly reciprocal, and
$\iota_{\rm rec}(\ell(\RR))=2$ otherwise.

Let $D^-$ and $D^+$ be nonempty properly immersed closed locally
convex subsets of $M$. For every $s> 0$, we denote by $\N_{D^-,\,D^+}
(s)$ the cardinality of the set of common perpendiculars from $D^-$ to
$D^+$ with length at most $s$, considered with multiplicities (see
\cite{ParPau17ETDS} or \cite{ParPau25Ingham} for precisions).  The
following result is the main tool in the proofs of this note.

\btheo\label{theo:1geoddivreal} Let $D^-$ be a Margulis cusp
neighbourhood in $M$ and let $D^+$ be the image of a divergent
geodesic in $M$.
Then as $s\ra+\infty$, we have
\[
\N_{D^-,D^+}(s)=
\frac{\Ga(\frac{n}{2})\,\iota_{\rm rec}(D^+)\,\Vol\partial D^-}
  {2\,\sqrt{\pi}\;\Ga(\frac{n+1}{2})\,m(D^+)\,\Vol M}\;
  s\, e^{(n-1) \,s} + \bigO(e^{(n-1) s})\,.
\]
\etheo

\dem Let $\|\sigma^+_{D^-}\|$ be the total mass of the outer skinning
measure\footnote{See Section 3 of \cite{ParPau17ETDS} for the
definition, that we won't need here.} $\sigma^+_{D^-}$ of $D^-$.  By
\cite[Prop.~20 (2)]{ParPau17ETDS}, we have $\|\sigma^+_{D^-}\|=
2^{n-1}\Vol\partial D^-$. By \cite[Thm.~6]{ParPau25Ingham}, we have
\[
\N_{D^-,D^+}(s)=
\frac{\Ga(\frac{n}{2})\,\iota_{\rm rec}(D^+)\,\|\sigma^+_{D^-}\|}
  {2^n\,\sqrt{\pi}\;\Ga(\frac{n+1}{2})\,m(D^+)\,\Vol M}\;
  s\, e^{(n-1) \,s} + \bigO(e^{(n-1) s})\,.
\]
Theorem \ref{theo:1geoddivreal} follows.
\cqfd

\medskip
The boundary at infinity of $\hnr$ is $\partial_\infty\hnr=\RR^{n-1}
\cup\{\infty\}$. For all distinct $x,y\in \partial_\infty\hnr$, let
$\wt\ell_{x,y}:\RR\ra\hnr$ be any geodesic line in $\hnr$ with points
at infinity $x= \wt\ell_{x,y} (-\infty)$ and $y=\wt\ell_{x,y}
(+\infty)$, unique up to translation at the source. For every discrete
group of isometries $\Ga$ of $\hnr$ with finite covolume such that $x$
and $y$ are parabolic fixed points of $\Ga$, we denote by $\ell_{x,y}=
\Ga\, \wt\ell_{x,y}$ the divergent geodesic in $\Ga\bs\hnr$ that is
the image of $\wt\ell_{x,y}$ under the quotient mapping.  Using the
terminology of \cite{McMullen21}, the image $\ell_{x,y}(\RR)$ of the
divergent geodesic $\ell_{x,y}$, endowed with its ``orientation'' from
$x$ to $y$, pushforward by $\ell_{x,y}$ of the orientation of $\RR$, is
a {\em degree $1$ modular symbol} in $\Ga\bs\hnr$, that we denote by
\[
\sym_\Ga(x,y)\;.
\]
Note that $\sym_\Ga(x,y)=\sym_\Ga(y,x)$ if and only if $\ell_{x,y}$ is
reciprocal in $\Ga\bs\hnr$, in which case $\sym_\Ga(x,y)$ will be
called a {\it reciprocal modular symbol}. We set
\[
\iota_{\Ga,\rm rec}(x,y)= \iota_{\rm rec}(\ell_{x,y}(\RR))\quad
     {\text and} \quad m_\Ga(x,y)=m(\ell_{x,y}(\RR))\,.
\]
Note that $\iota_{\Ga,\rm rec}(x,y)$ and $m_\Ga(x,y)$ are constant on
the $\Ga$-orbits of pairs $\{x,y\}$ of distinct parabolic fixed points
of $\Ga$.

\medskip
For $\KK\in\{\RR,\CC\}$, we denote by $\begin{bsmallmatrix}
a & b \\ c & d\end{bsmallmatrix}\in \PSL_2(\KK)$ the image of
$\begin{psmallmatrix}a & b \\ c & d\end{psmallmatrix}\in\SL_2(\KK)$.

\subsection{Counting modular symbols: the rational case}
\label{subsec:ratcase}
  
We identify as usual $\RR^2$ with $\CC$. Note that $\hdr=\big(\,
\{z\in\CC:\Im z>0\},\,ds^2=\frac{|dz|^2}{(\Im z)^2} \,\big)$. The
group $\PSLR$ acts on $\hdr\cup \partial_\infty \hdr$ by the
homographies $g\cdot z= \frac{az+b}{cz+d}$ for all $z\in \PP^1(\CC)$
and $g=\begin{bsmallmatrix} a & b \\ c & d\end{bsmallmatrix}\in
\PSL_2(\RR)$ with the usual convention when $z=\infty, -\frac{d}{c}$,
and acts faithfully by isometries on $\hdr$. The {\it modular group}
$\Ga_\QQ=\PSLZ$ is an arithmetic lattice in $\PSL_2(\RR)$. It acts
transitively on its set of parabolic fixed points $\PP^1(\QQ)=\QQ\cup
\{\infty\}$ in $\partial_\infty\hdr$. The stabiliser $\Ga_{\QQ,\infty}$
of $\infty$ in $\Ga_\QQ$ consists of the translations $z\mapsto z+k$
with $k\in\ZZ$.

\btheo\label{theo:rationalcount} Let $\Ga$ be a finite index subgroup
of $\Ga_\QQ=\PSLZ$, and let $\Ga_\infty$ be the stabiliser of $\infty$
in $\Ga$.  For all distinct $x,y\in\QQ\cup\{\infty\}$, as $\epsilon>0$
tends to $0$, we have
\[
\card\big(\Ga_{\infty}\bs\big\{ \{\alpha,\beta\}\in \Ga\cdot\{x,y\} :
\ |\beta-\alpha| \ge \epsilon\big\}\big) =
-\frac{6\;\iota_{\Ga,\rm rec}(x,y)\;[\Ga_{\QQ,\infty}:\Ga_\infty]}
{\pi^2\;[\Ga_{\QQ}:\Ga]}\;
\frac{\ln \epsilon}\epsilon+\bigO(\epsilon^{-1})\,.
\]
\etheo

For instance, for every $N\in\NN\ssm\{0,1\}$, with $\Ga=
\{\begin{bsmallmatrix} a & b \\ c & d\end{bsmallmatrix}\in \PSLZ:
c\equiv 0\!\!\mod N\}$ the Hecke congruence sugbroup of level $N$ of
$\Ga_\QQ=\PSLZ$, we have $[\Ga_{\QQ,\infty}:\Ga_\infty]=1$ and
$[\Ga_{\QQ}:\Ga]= N \prod_{p\mid N} \big(1+\frac{1}{p}\big)$ by
\cite[Prop.~1.43 (1)]{Shimura71} (as usual, the index $p$ ranges over
primes).

\medskip
\dem Note that $\Ga$ and $\Ga_\QQ$ have the same sets of parabolic
fixed points. We may hence apply Theorem \ref{theo:1geoddivreal} with
$M=\Ga\bs\hdr$, with $D^-=\Ga_{\infty}\bs B_\infty$ the image of
$B_\infty$ in $M$ (which is a Margulis cusp neighbourhood in $M$), and
with $D^+= \ell_{x,y}(\RR)=\Ga\, \wt\ell_{x,y}(\RR)$ (which is the
image of a divergent geodesic in $M$).  We have
\[
\vol M=[\Ga_{\QQ}:\Ga] \vol(\Ga_\QQ\bs \hdr)
=[\Ga_{\QQ}:\Ga]\;\frac\pi 3\, 
\]
and $\Vol\partial D^- = [\Ga_{\QQ,\infty}:\Ga_\infty] \Vol(
\Ga_{\QQ,\infty} \bs\partial B_\infty)=[\Ga_{\QQ,\infty}:\Ga_\infty]$.
Since the fixed point sets in $\hdr$ of the elliptic elements of
$\Ga_\QQ$ are singletons, the pointwise stabiliser in $\Ga_\QQ$ hence
in $\Ga$ of any geodesic line in $\hdr$ is trivial.  Therefore
$m(D^+)=1$.

If $\alpha,\beta\in\RR$ satisfy $|\alpha-\beta|<2$, then the length of
the common perpendicular from $B_\infty$ to $\wt\ell_{\alpha,\beta}
(\RR)$ is $\ln\big(\frac2{|\alpha-\beta|}\big)$. Since the stabilizer
in $\Ga_\QQ$ hence in $\Ga$ of a nontrivial geodesic segment is
trivial, the multiplicity of such a common perpendicular is $1$ (see
\cite[\S 3.3]{ParPau17ETDS}). Using these observations and Theorem
\ref{theo:1geoddivreal} with $n=2$, since $\Ga(\frac
32)=\frac{\sqrt\pi}2$, as $\epsilon>0$ tends to $0$, we obtain
\begin{align*}\label{eq:fareycountbydist}
  &\card\big(\Ga_\infty\bs\big\{ \{\alpha,\beta\}\in \Ga\cdot\{x,y\} : \
  |\beta-\alpha| \ge \epsilon\big\}\big) \nonumber\\
  =&\card\Big(\Ga_\infty\bs\Big\{ \{\alpha,\beta\}\in \Ga\cdot\{x,y\}  :
  d\big(B_\infty, \wt\ell_{\alpha,\beta}(\RR)\big)\le
  \ln\frac 2\epsilon\Big\}\Big)\nonumber\\
=&\,\N_{D^-,D^+}\Big(\ln\frac 2\epsilon\,\Big)\nonumber=
\frac{1\cdot \iota_{\Ga,\,{\rm rec}}(x,y)\cdot [\Ga_{\QQ,\infty}:\Ga_\infty]}
{2\sqrt\pi\,\frac{\sqrt\pi}2\cdot1\cdot[\Ga_{\QQ}:\Ga]\;\frac\pi 3}
\ln\Big(\frac 2\epsilon\Big) \frac 2\epsilon
+\bigO\Big(\frac 2\epsilon\Big)\,.
\end{align*}
Theorem \ref{theo:rationalcount} follows by simplification.
\cqfd

\subsection{Counting modular symbols: the imaginary quadratic case}
\label{subsec:imaquad}
  
The group $\PSL_2(\CC)$ acts on $\partial_\infty\htr=\PP^1(\CC)=
\CC\cup \{\infty\}$ by homographies (Möbius transformations) $g\cdot
z= \frac{az+b}{cz+d}$ for all $g=\begin{bsmallmatrix} a & b \\ c &
d\end{bsmallmatrix}\in \PSL_2(\CC)$ and $z\in\PP^1(\CC)$ with the
usual convention when $z=\infty, -\frac{d}{c}$. It acts faithfully
isometrically on
\[
\htr=\Big(\{(z,t)\in\CC\times \RR:t>0\}, \;\;ds^2=
\frac{|dz|^2+dt^2}{t^2}\Big)
\]
by Poincaré's extension.  

In this Subsection \ref{subsec:imaquad}, let $K$ be an imaginary
quadratic number field. Let $\OOO_K$, $D_K$, $h_K$, $\zeta_K$ be as in
the introduction. Recall that the group of units $\OOO_K^\times$ of
$\OOO_K$ is finite, and it is equal to $\{\pm 1\}$ unless $D_K=-4,-3$.
Let $\icg K$ be the group of ideal classes of $\OOO_K$, whose order is
$h_K$. For every ideal $\aaa$ of $\OOO_K$, recall that there exist
$a,b\in\OOO_K$ such that $\aaa=a\OOO_K+ b\OOO_K$, and we denote by
$[\aaa]$ the class of $\aaa$ in $\icg K$.  The identity element of
$\icg K$ is the principal class $[\OOO_K]$.

The {\it Bianchi group} $\Ga_K= \PSL_2(\OOO_K)$ is an arithmetic
lattice in $\PSL_2(\CC)$. The quotient space $M_K=\Ga_K\bs\htr$ is
hence a finite volume noncompact complete real hyperbolic
$3$-orbifold. The Bianchi group $\Ga_K$ acts with $h_K$ orbits on its
set of parabolic fixed points $\PP^1(K) =K\cup \{\infty\}$ in
$\partial_\infty\htr=\PP^1(\CC)$~: the map
\begin{equation}\label{cuspclass}
\Ga_K\cdot[x_0:x_1]\mapsto [x_0 \OOO_K + x_1\OOO_K]
\end{equation}
is a bijection\footnote{See for instance \cite[\S 7,
  Th.~2.4]{ElsGruMen98}} from $\Ga_K\bs\PP^1(K)$ to $\I_K$. The
stabiliser of $\infty$ in $\Ga_K$ is
\[
\Ga_{K,\infty}=\bigg\{\begin{bmatrix} a & b\\0&a^{-1}\end{bmatrix}\in
\Ga_K: a\in\OOO_K^\times\,,\; b\in\OOO_K\bigg\}\,.
\]
Note that $m_{\Ga_K}(x,y)$ is not constant when $x$ and $y$ vary among
the distinct elements of $K$. For example, we have $m_{\Ga_K}(0,
\infty) =\frac{|\OOO_K^\times|}2$, $m_{\Ga_K}(1,-1)=2$ and $m_{\Ga_K}
(\frac 13,\infty)=1$. As for $\iota_{\Ga_K,\,{\rm rec}} (x,y)$ (see
Subsection \ref{subsec:ford}), an explicit arithmetic value of
$m_{\Ga_K}(x,y)$ as $x$ and $y$ vary in $\PP^1(K)$ does not seem to be
available. See Examples (\hyperlink{exampfarid1}{1}) to
(\hyperlink{exampfarid3}{3}) in Section \ref{sec:cuspverymax} for
other examples of computation of $\iota_{\Ga_K,\,{\rm rec}} (x,y)$ and
$m_{\Ga_K}(x,y)$ for some $x,y\in\PP^1(K)$.

\btheo\label{theo:iqcount} Let $\Ga$ be a finite index subgroup of
$\Ga_K=\PSL_2(\OOO_K)$, and let $\Ga_\infty$ be the stabiliser of
$\infty$ in $\Ga$.  For all distinct $x,y\in K\cup\{\infty\}$, as
$\epsilon>0$ tends to $0$, we have
\begin{align*}
&\card\big(\Ga_{\infty}\bs\big\{ \{\alpha,\beta\}\in \Ga\cdot\{x,y\} :
|\beta-\alpha|\ge\epsilon\big\}\big)\\=&\;
\frac{4\;\pi\;\iota_{\Ga,\,{\rm rec}}(x,y)\;[\Ga_{K,\infty}:\Ga_\infty]}
{|\OOO_K^\times|\;D_K\;\zeta_K(2)\;m_\Ga(x,y)\;[\Ga_{K}:\Ga]}\;
\frac{\ln \epsilon}{\epsilon^2}+\bigO(\epsilon^{-2}).
\end{align*}
\etheo

\dem Again $\Ga$ and $\Ga_K$ have the same sets of parabolic fixed
points. As in the proof of Theorem \ref{theo:rationalcount}, we apply
Theorem \ref{theo:1geoddivreal} with $n=3$, with $M=\Ga\bs\htr$, with
$D^-=\Ga_{\infty}\bs B_\infty$, and with $D^+=\ell_{x,y}(\RR)$.  By
Humbert's volume formula,\footnote{See for instance \cite[\S 8.8 and
  \S 9.6]{ElsGruMen98}.}  we have
\[
\Vol(M)=[\Ga_K:\Ga]\,\Vol(\Ga_K\bs\htr)
=[\Ga_K:\Ga]\,\frac{|D_K|^{3/2}\,\zeta_K(2)}{4\,\pi}\;.
\]
The index $[\Ga_{K,\infty}:\OOO_K]$ in $\Ga_{K,\infty}$ of its
unipotent subgroup consisting in the translations by elements of
$\OOO_K$ is equal to $\frac{|\OOO_K^\times|}{2}$. By for instance the
area computation in the proof of \cite[Lemma 6]{ParPau11BLMS}, since
$\OOO_K$ is generated as a $\ZZ$-module by $1$ and $\frac{D_K +
  i\sqrt{|D_K |}}2$, we have
\begin{align*}
\Vol(\partial D^-)&=[\Ga_{K,\infty}:\Ga_\infty]
\Vol(\Ga_{K,\infty}\bs \partial B_\infty)=[\Ga_{K,\infty}:\Ga_\infty]
  \frac{2}{|\OOO_K^\times|}\Vol(\OOO_K\bs \CC)
 \\&=[\Ga_{K,\infty}:\Ga_\infty]\;\frac{\sqrt{|D_K|}}{|\OOO_K^\times|}\,.
\end{align*}
Only finitely many $\Ga_{K,\infty}$-orbits (hence $\Ga_\infty$- orbits)
of geodesic lines $\wt\ell_{\infty,z}(\RR)$ for $z\in\CC$ have
nontrivial stabilizers in $\Ga_K$ (hence in $\Ga$). A given geodesic
line $\wt\ell_{\infty,z}(\RR)$ meets perpendicularly at most
$\bigO(t)$ elements of $\Ga\,\wt\ell_{x,y}(\RR)$ at distance at most
$t$ from $B_\infty$.  Hence only linearly many $\Ga_\infty$-orbits of
common perpendiculars between $B_\infty$ and elements of $\Ga\,\wt
\ell_{x,y}(\RR)$ have multiplicity different from $1$. As in the proof
of Theorem \ref{theo:rationalcount}, when $\epsilon>0$ tends to $0$,
we have
\begin{align*}
  & \card\big(\Ga_\infty\bs\big\{ \{\alpha,\beta\}\in \Ga\cdot\{x,y\} :  \
  |\beta-\alpha| \ge \epsilon\big\}\big)=
  \N_{D^-,D^+}\Big(\ln\frac 2\epsilon\Big)+
  \bigO\Big(\ln\frac 2\epsilon\Big)\\=&\,
  \frac{\Ga(\frac{3}{2})\;\iota_{\Ga,\,{\rm rec}}(x,y)\;
    [\Ga_{K,\infty}:\Ga_\infty]\;\frac{\sqrt{|D_K|}}
    {|\OOO_K^\times|}}{2\,\sqrt{\pi}\;\Ga(2)\;m_\Ga(x,y)\;
    [\Ga_K:\Ga]\,\frac{|D_K|^{3/2}\;\zeta_K(2)}
    {4\,\pi}}\; \Big(\ln\frac 2\epsilon\Big)
  \Big(\frac 2\epsilon\Big)^2  + \bigO(\epsilon^{-2})\,.
\end{align*}
Theorem \ref{theo:iqcount} follows by simplification, using the values
$\Ga(2)=1$ and $\Ga(\frac 32)=\frac{\sqrt\pi}2$.  \cqfd

\subsection{Proofs of Theorem \ref{theo:realfareycount} and
  Corollary \ref{cor:farey}}
\label{subsec:proofintro}
  
In this Subsection, let $K$ be as in the introduction, either $\QQ$ or
an imaginary quadratic  number field. We are now ready to prove Theorem
\ref{theo:realfareycount}, by restricting Theorems
\ref{theo:rationalcount} and \ref{theo:iqcount} to the case when
$x,y\in \PP^1(K)$ are Farey neighbours.  Since $\Ga_K=\PSL_2(\OOO_K)$
has infinitely many orbits on the set of pairs of points of
$\PP^1(K)$, Lemma \ref{lem:fromarithtogeom} below implies in
particular that there are lots of unordered pairs of distinct elements
$x,y\in K\cup\{\infty\}$ that are not Farey neighbours, though
Theorems \ref{theo:rationalcount} and \ref{theo:iqcount} still
apply. See Examples (\hyperlink{exampfarid1}{1}) to
(\hyperlink{exampfarid3}{3}) in Section \ref{sec:cuspverymax} for
explicit examples of degree $1$ modular symbols, that are not pairs
of Farey neighbours, in the imaginary quadratic case. The key
translation between the arithmetics and the geometry is the following
elementary lemma (that will be more involved in the quaternionic
case).

\blemm\label{lem:fromarithtogeom} Two distinct elements $\alpha,\beta
\in \PP^1(K)=K\cup \{\infty\}$ are Farey neighbours if and only if
there exists $\ga\in\Ga_K=\PSL_2(\OOO_K)$ such that $\ga\cdot \infty=
\alpha$ and $\ga\cdot 0= \beta$.
\elemm

\dem Let $\alpha,\beta\in \PP^1(K)$ be distinct. If there exists
$\ga=\begin{bsmallmatrix} p&r\\q&s\end{bsmallmatrix} \in\Ga_K$ such
that $\ga\cdot \infty= \alpha$ and $\ga\cdot 0= \beta$, then $p,q,r,s
\in\OOO_K$, $\alpha=\frac pq$, $\beta=\frac rs$ and $ps-qr=1$, hence
$\alpha,\beta$ are Farey neighbours.

Conversely, if $\alpha,\beta$ are Farey neighbours, let $p,q,r,s\in
\OOO_K$ be such that $\alpha=\frac pq$, $\beta=\frac rs$ and
$|\,ps-qr\,|=1$. Then $u=ps-qr\in\OOO_K^\times$. Hence if $\ga=
\begin{bsmallmatrix} pu^{-1}&r\\qu^{-1}&s\end{bsmallmatrix}$,
then $\ga\in\Ga_K$ and we have $\ga\cdot \infty= \alpha$ and
$\ga\cdot 0= \beta$.
\cqfd

\brema\label{rem:orbitalequivFareyneigh}{\rm
\hypertarget{orbitalequivFareyneighi}{(i)} Lemma
\ref{lem:fromarithtogeom}, and the fact that the ideal classes
associated with $\infty =\frac 10$ and $0=\frac 01$ by Equation
\eqref{cuspclass} are the principal one $[\OOO_K]$, imply that if
$x=[x_0:x_1] \in \PP^1(K)$ and $y=[y_0:y_1]\in \PP^1(K)$ are Farey
neighbours, then their two associated ideal classes
$[x_0\OOO_K+x_1\OOO_K]$ and $[y_0\OOO_K + y_1\OOO_K]$ are the
principal one.

\medskip \hypertarget{orbitalequivFareyneighii}{(ii)} In this remark,
assume that $\OOO_K$ is Euclidean, that is, that $D_K=1$ ($K=\QQ$) or
$D_K=-3,-4,-7,-8,-11$. Let $n_K=[K:\QQ]$, let $\Ga$ be a finite index
subgroup of $\Ga_K$ and let $X_\Ga=\Ga\,\bs(\HH^{n_k}_\RR\cup K\cup \{
\infty\})$ be the cuspidal compactification of $\Ga\bs\HH^{n_k} _\RR$.
By \cite[\S 1.2]{Manin72} when $K=\QQ$ and as extended in \cite[\S
  2.2]{Cremona84} otherwise, the degree $1$ modular symbols
$\sym_\Ga(x,y)$ for distinct $x,y\in K\cup\{\infty\}$, with the cusp
points from which they start (respectively end) added at their
beginning (respectively end), are integral $1$-cycles when $\Ga x=\Ga
y$, and define real $1$-cycles when $\Ga x\neq\Ga y$, whose homology
classes in $H_1(X_\Ga,\RR)$ we denote by $[\sym_\Ga(x,y)]$. By Lemma
\ref{lem:fromarithtogeom}, precisely when $x$ and $y$ are Farey
neighbours, the homology classes $[\sym_\Ga (x,y)]$ are called {\it
  distinguished classes} in \cite[\S 1.5] {Manin72} when $K=\QQ$ and
{\it special classes} in \cite[\S 2.2] {Cremona84} otherwise. Since
$\OOO_K$ is Euclidean, these finitely many classes generate the real
vector space $H_1(X_\Ga,\RR)$ by \cite[Prop.~1.6 a)]{Manin72} when
$K=\QQ$ and \cite[page 287, lines - 9 to - 5]{Cremona84} otherwise.

When $\OOO_K$ is not Euclidean, the same arguments prove that the
homology classes of the degree $1$ modular symbols corresponding to
the $1$-edges of the dual ideal tesselation of the Ford-Voronoi
tesselation of $\htr$ (whose $2$-skeleton is the Mendoza
$\Ga_K$-invariant spine of $\htr$, see \cite[\S 4]{BesSav12} and
Section \ref{sec:cuspverymax}) generate the real vector space
$H_1(X_\Ga,\RR)$.  } \erema

\noindent{\bf Proof of Theorem \ref{theo:realfareycount}.}  In order
to prove Claim (\hyperlink{mainintro1}{1}), we apply Theorem
\ref{theo:rationalcount} to $\Ga=\Ga_\QQ$ and $(x,y)=(0,\infty)$, so
that by Lemma \ref{lem:fromarithtogeom}, we have $\fn_{\QQ}=
\Ga\cdot\{x,y\}$.  The divergent geodesic $\ell_{0,\infty}$ in
$\Ga_\QQ\bs\hdr$ is reciprocal since the elliptic element
$\begin{bsmallmatrix} 0&-1\\1&\ 0\end{bsmallmatrix}\in\Ga_\QQ$ of
order $2$ exchanges the two points at infinity $0$ and $\infty$ of the
geodesic line $\wt\ell_{0,\infty}(\RR)$. Hence $\iota_{\Ga_\QQ,\rm
  rec} (0, \infty)=1$.  Two images of the geodesic line
$\wt\ell_{0,\infty} (\RR)$ under $\Ga_\QQ$ either coincide or do not
intersect in $\hdr$. If the two endpoints of such an image are
different from $\infty$, then their distance is at most one. Thus,
except the images of $\wt\ell_{0,\infty}(\RR)$ under
$\Ga_{\QQ,\infty}$, every image of $\wt\ell_{0,\infty}(\RR)$ under an
element of $\Ga$ has one and only one translate modulo
$\Ga_{\QQ,\infty}= \ZZ$ both of whose points at infinity lie in the
unit interval $[0,1]$.  Thus Claim (\hyperlink{mainintro1}{1}) of
Theorem \ref{theo:realfareycount} follows from Theorem
\ref{theo:rationalcount}.

\medskip
In order to prove Claim (\hyperlink{mainintro2}{2}) of
Theorem \ref{theo:realfareycount}, we apply similarly Theorem
\ref{theo:iqcount} with $\Ga=\Ga_K$ and $(x,y)=(0,\infty)$, so that by
Lemma \ref{lem:fromarithtogeom}, we have $\fn_{K}= \Ga\cdot\{x,y\}$.  Note
that the pointwise stabiliser in $\Ga_K$ of the geodesic line
$\wt\ell_{0,\infty}(\RR)$ has cardinality $\frac{|\OOO_K^\times|}2$,
hence $m_{\Ga_K}(0,\infty)=\frac{|\OOO_K^\times|}2$. Note that the
divergent geodesic $\ell_{0,\infty}$ is reciprocal as in the
rational case, hence $\iota_{\Ga_K,\rm rec}(0,\infty)=1$. The index in
$\Ga_{K,\infty}$ of its unipotent subgroup of translations by $\OOO_K$
is equal to $\frac{|\OOO_K^\times|}2$. Hence replacing the quotient
modulo $\Ga_{K,\infty}$ in the left-hand side of the formula in
Theorem \ref{theo:iqcount} by the quotient modulo $\OOO_K$ amounts to
multiplying the right-hand side by $\frac{|\OOO_K^\times|}2$. \cqfd

\bigskip
\noindent{\bf Proof of Corollary \ref{cor:farey}.} Note that Equation
\eqref{eq:defiFareyNeigh} is equivalent when $qs\neq 0$, to the
equation $\big|\frac pq-\frac rs\big|=\frac 1{|qs|}$. For every
$N\in\NN\ssm\{0\}$, the map $(p,q,r,s)\mapsto \{\frac pq,\frac rs\}$
from the set
\[
\big\{(p,q,r,s)\in\ZZ^4: ps-qr=1\,,\ 0\le p\le
q\,,\ 0\le r\le s\,,\ 0<qs\le N\big\}
\]
to the set $\big\{\{\alpha,\beta\} \in \fn_{\QQ}: \alpha,\beta
\in\interval 01\,,\ |\alpha-\beta|\ge \frac 1N\big\}$ is easily
checked to be $2$-to-$1$. Hence Corollary \ref{cor:farey} follows from
Theorem \ref{theo:realfareycount} (\hyperlink{mainintro1}{1}) with
$\epsilon=\frac{1}{N}$.
\cqfd

\medskip The next pictures illustrate Theorem
\ref{theo:realfareycount} (\hyperlink{mainintro1}{1}) and Corollary
\ref{cor:farey}. The blue curve on the left represents the graph of
the map $N\mapsto \card\big\{ \{\alpha,\beta\}\in \fn_{\QQ} : \alpha,
\beta\in\interval 01\,,\ |\beta-\alpha| \ge \frac1N\big\}$.  The
orange one represents the graph of the map $N\mapsto \frac6{\pi^2}
\,N\,\ln N$. Note that the two graphs diverge slowly one from the other
since there is only a logarithmic factor between the main term and the
error term. The picture on the right represents the ratio map, slowly
converging to $1$.

\begin{center}
  \includegraphics[width=7cm]{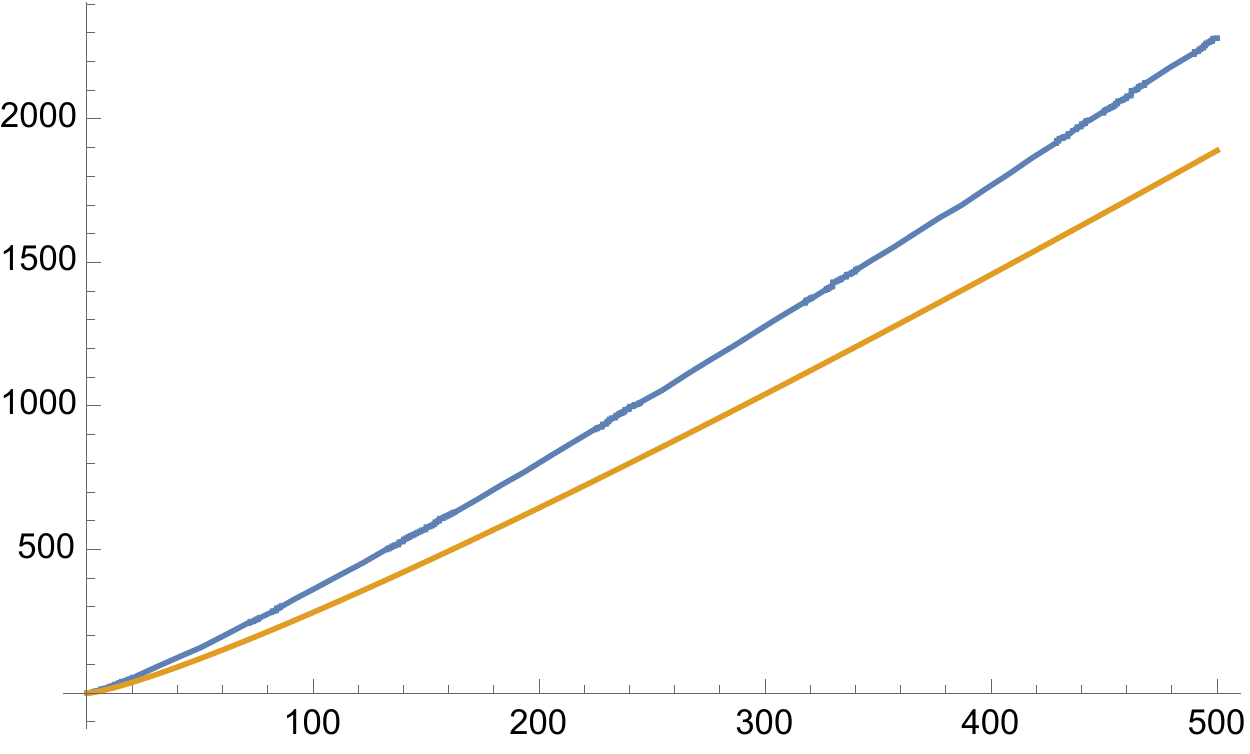}\quad
  \includegraphics[width=7cm]{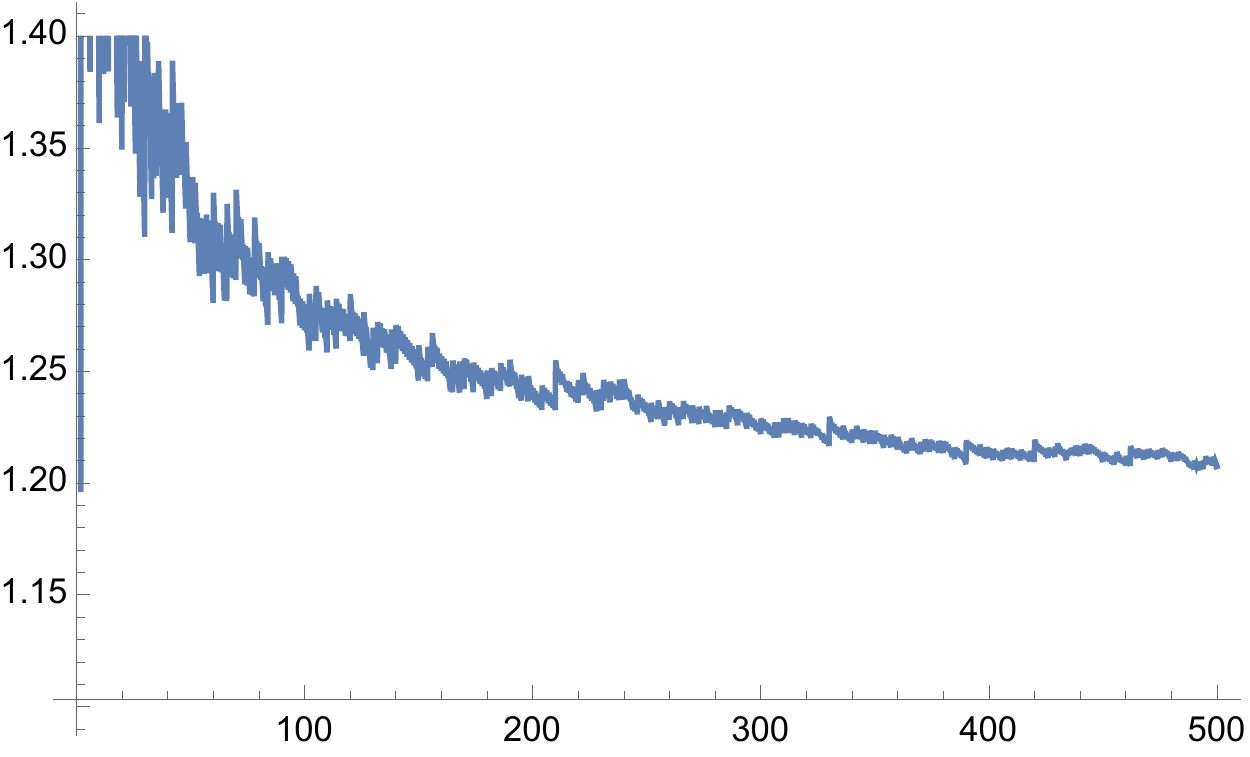}
\end{center}

\subsection{Ford circles, Farey neighbours and modular symbols}
\label{subsec:ford}

Assume in this Subsection that $K=\QQ$. Being Farey neighbours in
$\PP^1(\QQ)$ has a well-known geometric characterisation, that we
recall as a motivation for Section \ref{sec:cuspverymax}. For every
$\alpha\in \PP^1(\QQ)$, if $\alpha=\frac pq$ with $p,q\in\ZZ$
relatively prime and $q> 0$, let $B_\alpha$ be the intersection with
$\hdr$ of the closed Euclidean ball of center $x+\frac i{2q^2}$ and
radius $\frac 1{2q^2}$. The boundary of this disc is called the {\em
  Ford circle}\footnote{Ford himself calls them {\it Speiser
  circles.}} of $\frac pq\in \QQ$, see the picture below. Let
$B_\infty=\{z\in\hdr:\Im z\ge 1\}$.  The family
$(B_\alpha)_{\alpha\in\PP^1(\QQ)}$ is the unique
$\PSL_2(\ZZ)$-equivariant family of maximal horoballs with pairwise
disjoint interiors. Two distinct $\alpha,\beta\in\PP^1(\QQ)$ are Farey
neighbours if and only if the horoballs $B_\alpha$ and $B_\beta$ are
tangent, or if and only if their Ford circles are tangent, see for
instance \cite[page 12]{Zullig28} that was published before
\cite{Ford38}.

\begin{center}
\includegraphics[width=12cm]{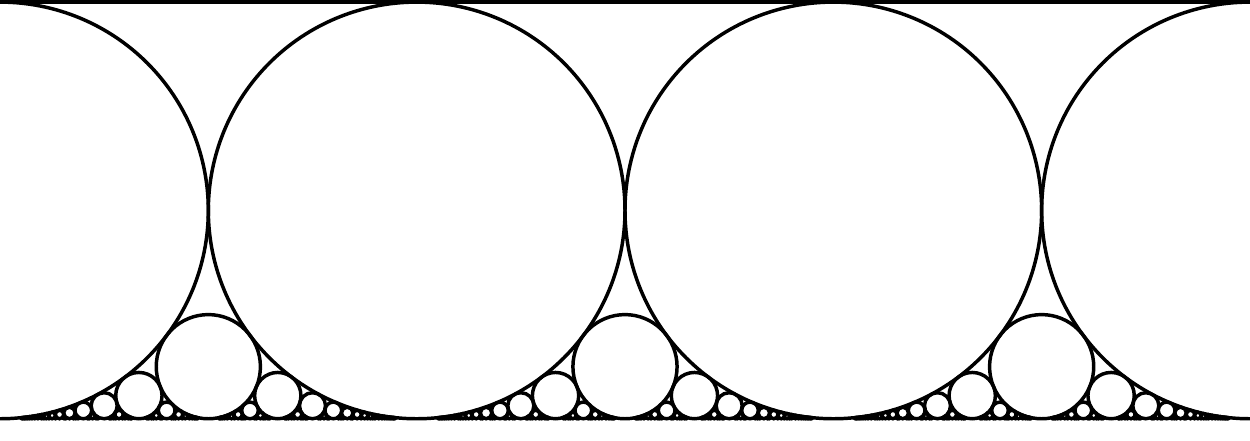}
\end{center}

\brema\label{rem:ford}{\rm
\hypertarget{fordi}{(i)} {\bf Another counting of degree $1$ modular
  symbols. } For all distinct $x,y\in\PP^1(\QQ)$, the hyperbolic
distance $d(B_x,B_y)$ is a natural complexity for the degree $1$
modular symbol $\sym_{\Ga_\QQ}(x,y)$, see \cite{ParPauSay25}.  For
every $T>0$, let
\[
\Sym(T)=\{\sym_{\Ga_\QQ}(x,y):
x,y\in\PP^1(\QQ),\; x\neq y,\; d(B_x,B_y)\le T\}
\]
and let $\Symrec(T)$ be the subset of $\Sym(T)$ that consists of its
reciprocal modular symbols. Let $D^-=D^+$ be the Margulis
neighbourhood of the cusp of $\Ga_\QQ\bs\hdr$ defined as the image of
any $B_\alpha$ for $\alpha\in\PP^1(\QQ)$ under the quotient mapping
$\hdr\ra\Ga_\QQ\bs\hdr$. Its hyperbolic area is $1$. By
\cite[Cor.~21]{ParPau17ETDS}, there exists $\kappa>0$ such that as
$T\ra+\infty$, we have
\begin{align*}
\card\Sym(T)&=\N_{D^-,\,D^+}(T)=\frac{2^{2-1}(2-1)\Vol(D^-)\Vol(D^+)}
{\Vol(\SSS^{2-1})\Vol(\Ga_\QQ\bs\hdr)}\, e^T(1+e^{-\kappa T})
\\&=\frac 3{\pi^2}\,e^T(1+e^{-\kappa T})\,.
\end{align*}

For all distinct $x,y\in\PP^1(\QQ)$, the degree $1$ modular symbol
$\sym_{\Ga_\QQ}(x,y)$ for $\Ga_\QQ$ is reciprocal if and only if the
geodesic line $\wt\ell_{x,y}(\RR)$ in $\hdr$ intersects the orbit
$\Ga_\QQ \cdot i$. This (unique) point of intersection is the midpoint
of the common perpendicular of $B_x$ and $B_y$. Thus, for every $T>0$,
the number of reciprocal modular symbols of complexity at most $T$
equals the number of common perpendiculars in $\Ga_\QQ\bs \hdr$ from
$D^-$ as above to ${D'}^+=\Ga_\QQ\cdot i$ of length at most $\frac
T2$.  Since the stabilizer of $i$ in $\Ga_\QQ$ has cardinality $2$, by
\cite[Cor.~21]{ParPau17ETDS}, there exists $\kappa>0$ such that as
$T\ra+\infty$, we have
\begin{equation}\label{eq:symrecas}
\card\Symrec(T)=\N_{D^-,\,{D'}^+}\Big(\frac T2\Big)=
\frac{\Vol(D^-)}{2\Vol(\Ga_\QQ\bs\hdr)}\, e^{\frac T2}
(1+e^{-\kappa \frac T2})=\frac 3{2\pi}\,e^{\frac T2}(1+e^{-\kappa\frac T2})\,.
\end{equation}
Thus, the proportion of reciprocal modular symbols in $\Sym(T)$ is
equivalent to $2\pi\,e^{-\frac T2}$ as $T\ra+\infty$.

\medskip
\hypertarget{fordii}{(ii)} {\bf Relationship between counting modular
  symbols and the primitive circle problem. } For
every $n\in\NN$, a {\it representation by primitive sums of two
  squares} of $n$ is a pair $(p,q)\in\ZZ^2$ with $p,q$ coprime such
that $n=p^2+q^2$.  Let us denote by $r_{\rm prim }(n)\,$ their
number. As $N\ra+\infty$, we have
\begin{equation}\label{eq:averprimrep}
\sum_{n=1}^N r_{\rm prim}(n)=\frac{6}{\pi}\,N+\bigO\big(\sqrt{N}\,\big)\,,
\end{equation}
see \cite[Eq.~(1.1)]{Wu02} (and Theorem 1 in loc.~cit.~for a better
error term conditionally to the RH). Let us prove that Equation
\eqref{eq:symrecas} follows from Equation \eqref{eq:averprimrep}.

For every $\ga=\begin{bsmallmatrix} a & b \\ c & d\end{bsmallmatrix}
\in\Ga_\QQ$, the integers $c$ and $d$ are coprime, and we have
\begin{equation}\label{eq:ReImOrbi}
  \Re(\ga\cdot i)=\frac{ac+bd}{c^2+d^2}\quad\text{and}\quad
  \Im(\ga\cdot i)=\frac{1}{c^2+d^2}\,.
\end{equation}
In particular, the imaginary part of any element of $\Ga_\QQ\cdot i$
has the form $\frac{1}{n}$ for some $n\in\NN\smallsetminus\{0\}$ which
is a primitive sum of two squares.  Fixing coprime integers $c,d\in
\ZZ$ and a solution $(a,b)\in\ZZ^2$ of the gcd equation $ad-bc=1$, all
other solutions are $(a+kc,b+kd)$ with $k\in \ZZ$ by the uniqueness
property of the Bézout identity. Since $\frac{(a+kc)c+(b+kd)d}
{c^2+d^2} =\frac{ac+bd}{c^2+d^2}+k$, there exists a unique such
solution $(a',b')\in\ZZ^2$ such that $R(c,d)= \frac{a'c+b'd}{c^2+d^2}$
belongs to $[0,1[\,$. We define $\ga_{c,d} =\begin{bsmallmatrix}
\,a'&\,b' \\ c & d\end{bsmallmatrix}$. For every $n\in\NN$, given a
representation $(c,d)$ of $n$ by sums of two squares with $d\neq \pm
c$, there are $8$ representations of $n$ obtained by changing the order
and the signs of $c$ and $d$. Among these $8$ representations, the $4$
pairs $(c,d)$, $(-c,-d)$, $(d,-c)$, $(-d,c)$ do not change $R(c,d)$
(we have $\ga_{-c,-d}=\ga_{c,d}$ and $\ga_{d,-c}=\ga_{-d,c} =
\ga_{c,d} \circ \iota$ with $\iota= \begin{bsmallmatrix} 0 & -1 \\ 1 &
0\end{bsmallmatrix}$ fixing $i$). The $4$ remaining pairs obtained by
exchanging $c$ and $d$ change $R(c,d)$ into $1-R(c,d)$. The number of
coprime pairs $(c,d)\in\ZZ^2$ such that $R(c,d)\in\{0,\frac{1}{2},1\}$
is finite, since a vertical geodesic line in $\hdr$ meets at most one
point of $\Ga_\QQ\cdot i$. The number of all representations of
integers by primitive sums of squares of two equal or opposite
integers is finite (equal to $4$, since $c\in\NN$ is coprime to $\pm
c$ if and only if $c=\pm 1$). Thus as $T\ra+\infty$, by the standard
computation of the hyperbolic distance of a point of $\hdr$ to the
horizontal horosphere $\partial D^-=\{z\in\hdr: \Im z= 1\}$ and by
Equation \eqref{eq:averprimrep}, we have
\begin{align*}
  \N_{D^-,\,{D'}^+}\Big(\frac T2\Big)&= \card\big\{z\in\Ga_\QQ\cdot i: 0\le\Re
  z<1,\; \Im z\in [e^{- \frac{T}{2}},1]\;\big\}\\&
  =\frac{1}{4}\;\sum_{n=1}^{\big\lfloor e^{\frac{T}{2}}\big\rfloor} r_{\rm prim}(n)
  +\bigO\big(1\big)=\frac{3}{2\pi}\,e^{\frac{T}{2}}+
  \bigO\big(e^{\frac{T}{4}}\,\big)\,.
\end{align*}
This implies Equation \eqref{eq:symrecas}, as wanted, with an explicit value $\kappa=\frac 12$.

\medskip
\hypertarget{fordiii}{(iii)} {\bf On the computation of the
  reciprocity indexes. } Given a finite index subgroup $\Ga$ of
$\Ga_\QQ$, as $x$ and $y$ vary in $\PP^1(\QQ)$, finding an explicit
arithmetic value of the reciprocity index $\iota_{\Ga,\rm rec}(x,y)$
is somewhat delicate, even when $\Ga=\Ga_\QQ$. This also turns out to
be related to problems of representations by primitive sums of two
squares, as we now indicate.

By the diagonal $\Ga_\QQ$-invariance and by the transitivity of the
action of $\Ga_\QQ=\PSLZ$ on $\PP^1(\QQ)$, we only need to compute
$\iota_{\Ga_\QQ,\rm rec}(\infty,x)$ for $x\in K\cap[0,1[\,$. We have
$\iota_{\Ga_\QQ,\rm rec}(\infty,x)=1$ if and only if the geodesic line
$\wt \ell_{\infty,x}(\RR)$ meets the $\Ga_\QQ$-orbit of $i$, that is,
if and only if there exists $\ga\in\Ga_\QQ$ such that $\Re(\ga\cdot i)
=x$. Let us write $x=\frac pq$ with $p,q\in\ZZ$ coprime and $q>0$.
Note that for all $a,b,c,d\in \ZZ$ such that $ad-bc=1$, we have
$c^2+d^2>0$ and $(a^2+b^2)(c^2+d^2)-(ac+bd)^2=(ad-bc)^2=1$ by the
Diophantus identity.  Hence $ac+bd$ and $c^2+d^2$ are coprime by the
Bézout identity, so that by the uniqueness property of reduced
fractions, if $\frac pq=\frac{ac+bd}{c^2+d^2}$, then $q=c^2+d^2$ and
$p=ac+bd$.  By Equation \eqref{eq:ReImOrbi} and the discussion of
Claim (\hyperlink{fordii}{ii}), we hence have $\iota_{\Ga_\QQ,\rm rec}
(\infty, x) =1$ if and only if $q=c^2+d^2$ is a primitive sum of two
squares and $p=qR(c,d)$.

} \erema

\section{Bianchi cusps are very maximal}
\label{sec:cuspverymax}

Let $K$, $\OOO_K$, $\icg K$, $h_K$, $D_K$, $\Ga_K=\PSL_2(\OOO_K)$ and
$M_K=\Ga_K \bs \htr$ be as in Subsection \ref{subsec:imaquad}.  Let
$\radix$ be a square-free negative integer such that $K=\QQ(\radix)$,
with $D_K=4\radix$ if $\radix\equiv 2,3\bmod 4$ and $D_K=\radix$
otherwise.  When $h_K\neq 1$, there are elements of $\PP^1(K)$ that have no 
Farey neighbour as defined in Equation \eqref{eq:defiFareyNeigh}, by
Remark \ref{rem:orbitalequivFareyneigh}
(\hyperlink{orbitalequivFareyneighi}{i}). The aim of this Section is
to advertise a more general notion of Farey neighbours and to prove
that it solves this issue. We refer for instance to \cite{Mendoza80},
\cite[Chap.~7]{ElsGruMen98}, \cite[Sect. 4]{BesSav12} for background
material on this Section.

Two distinct $x,y\in\PP^1(K)$ are said to be {\it \genfar neighbours}
if for any $a,b,c,d\in \OOO_K$ such that $x=\frac ab$ and $y=\frac
cd$, we have
\begin{equation}\label{eq:fareyideal}
(a\OOO_K+b\OOO_K)(c\OOO_K+d\OOO_K)=(ad-bc)\OOO_K\,.
\end{equation}
This does not depend on the choices of $a,b,c,d$.  We refer to
Examples \ref{ex:rami} (\hyperlink{exampfarid1}{1}) to
(\hyperlink{exampfarid3}{3}) below for examples of \genfar
neighbours. After some remarks, we will recall the geometric
interpretation of this property, due to \cite{BesSav12}, that proves
without computation that being \genfar neighbours is a property
invariant by the diagonal action of $\Ga_K$ on the set of unordered
pairs $\{x,y\}$ of distinct elements of $\PP^1(K)$.

\brema\label{rem:remdefgenfar} {\rm
\hypertarget{remdefgenfar1}{(1)} If $x,y\in\PP^1(K)$ are Farey
neighbours, then they are \genfar neighbours. Indeed, $(a,b,c,d)=
(1,0,0,1)$ is a solution of Equation \eqref{eq:fareyideal}, hence
$\infty=\frac 10, 0=\frac 01$ are \genfar neighbours. Being \genfar
neighbours is invariant by $\Ga_K$, hence this Claim
(\hyperlink{remdefgenfar1}{1}) follows from Lemma
\ref{lem:fromarithtogeom}. For example, Equation \eqref{eq:fareyideal}
implies that the \genfar neighbours of $x=\infty =\frac 10$ are the
points $\frac cd$ with $c\OOO_K+ d\OOO_K=d\OOO_K$ or equivalently
$d\mid c$, hence are the points in $ \OOO_K$, that is, are its Farey
neighbours. By $\Ga_K$-invariance, the \genfar neighbours of an
element $x\in\PP^1(K)$ whose associated ideal class by Equation
\eqref{cuspclass} is principal are its Farey neighbours.  In
particular if there exist pairs of \genfar neighbours that are not
pairs of Farey neighbours, then $h_K\geq 2$.\footnote{See the comment
after Theorem \ref{theo:bianverymax} for the converse.}

\medskip \hypertarget{remdefgenfar2}{(2)} If $x,y\in\PP^1(K)$ are
\genfar neighbours, then their associated ideal classes are inverse
one of the other in the group $\icg K$ : by Equation
\eqref{eq:fareyideal}, if $a,b,c,d\in\OOO_K$ are such that $x=\frac
ab$ and $y=\frac cd$ are \genfar neighbours, then
$[a\OOO_K+b\OOO_K]^{-1} = [c\OOO_K+d\OOO_K]$. Hence if furthermore the
divergent geodesic $\ell_{x,y}$ is reciprocal, then $x$ and $y$ in
particular are in the same $\Ga_K$-orbit, thus have same associated
ideal class, which is either trivial or has order $2$ in the group
$\icg K$. We refer to Examples (\hyperlink{exampfarid1}{1}) to
(\hyperlink{exampfarid3}{3}) below for examples of such order $2$
ideal classes.

\medskip \hypertarget{remdefgenfar3}{(3)} Let $\Nr(\aaa)=
[\OOO_K:\aaa]$ be the norm of a nonzero ideal $\aaa$ of $\OOO_K$,
extended by multiplicativity to the norm of fractional ideals. For
every $a\in K$, let $\Nr(a)=\Nr(a\OOO_K)$.  Equation
\eqref{eq:fareyideal} implies that $\Nr(\OOO_K+x\OOO_K)
\Nr(\OOO_K+y\OOO_K)=\Nr(x-y)$, which is the equality case in the
inequalities with $c_1=c_2=1$ pages 10 and 11 of \cite{Mendoza80}.
}\erema

Let us turn to a geometric characterisation of being \genfar
neighbours in the complex case, that is analogous to the tangency
property of Ford circles discussed in Subsection \ref{subsec:ford}.
For every $\alpha\in\PP^1(K) \ssm\{\infty\}$, writing $\alpha=\frac
ab$ with any $a,b\in\OOO_K$, the {\em canonical (closed) horoball}
$B_{\alpha}$ in $\htr$ centered at $\alpha$ is the intersection with
$\htr$ of the Euclidean closed ball with center $\big(\alpha,
\frac{\Nr(a\OOO_K+ b\OOO_K)}{2\Nr(b)} \big)\in\htr$ and radius
$\frac{\Nr(a\OOO_K+b\OOO_K)}{2\Nr(b)}$. This does not depend on the
choices of $a$ and $b$.  Furthermore, the {\em canonical horoball} in
$\htr$ centered at $\infty$ is the already defined horoball $B_\infty=
\{(z,t) \in\htr:t\geq 1\}$.

This family $(B_x)_{x\in \PP^1(K)}$, constructed and studied in
\cite{Mendoza80}, is a $\Ga_K$-equivariant family of horoballs with
pairwise disjoint interiors.  In particular, for every $\ga\in\Ga_K$,
two canonical horoballs $B_x$ and $B_y$ with distinct $x,y\in
\PP^1(K)$ touch at one point (or equivalently are tangent) if and only
if $B_{\ga\cdot x}$ and $B_{\ga\cdot y}$ are tangent. The image of
$\bigcup_{x\in \PP^1(K)} B_x$ in the quotient orbifold
$M_K=\Ga_K\bs\htr$ is the union of closed Margulis cusp neighbourhoods
with pairwise disjoint interiors of the $h_K$ ends of $M_K$.

The geometric characterisation alluded to above, proving the
$\Ga_K$-invariance of being \genfar neighbours, is the following one, see \cite[Prop.~4.1]{BesSav12} for the proof.

\bprop\label{prop:BesSav}  Two distinct
elements $x,y\in \PP^1(K)$ are \genfar  neighbours if and only if the
canonical horoballs $B_x$ and $B_y$ are tangent. \cqfd
\eprop

\brema{\rm \hypertarget{remdefcanhoro1}{(1)} It follows from
Proposition \ref{prop:BesSav} and Remark \ref{rem:remdefgenfar}
(\hyperlink{remdefgenfar2}{2}) that for all distinct $x,y\in
\PP^1(K)$, if the canonical horoballs $B_x$ and $B_y$ are tangent,
then the ideal classes associated with $x$ and $y$ are inverse one of
the other.

\medskip
\hypertarget{remdefcanhoro2}{(2)} Since $\Nr(a\OOO_K+b\OOO_K)=
\Nr(b\OOO_K)$ implies that $a\OOO_K+b\OOO_K=b\OOO_K$ hence that $b\mid
a$ for all $a,b\in\OOO_K$, it follows from their construction that the
canonical horoballs that are tangent (and distinct) to the canonical
horoball $B_\infty$ are the ones centered at $c=\frac c1\in\OOO_K$,
confirming the example claim of Remark \ref{rem:remdefgenfar}
(\hyperlink{remdefgenfar1}{1}), by Proposition \ref{prop:BesSav}.  }
\erema

The main result of this Section \ref{sec:cuspverymax}, proving the
maximality of $(B_x)_{x\in \PP^1(K)}$ at all cusps, is the following
one.

\btheo\label{theo:bianverymax} Every element of
$\PP^1(K)$ has infinitely many \genfar  neighbours.
\etheo

By Proposition \ref{prop:BesSav}, we have an equivalent, more geometric
formulation of Theorem \ref{theo:bianverymax}.

\btheo For every $x\in \PP^1(K)$, the canonical horoball
$B_x$ is tangent to infinitely many elements of Mendoza's canonical
family $(B_x)_{x\in \PP^1(K)}$ of horoballs. \qed
\etheo

If $h_K\geq 2$, then any element of $\PP^1(K)$ whose associated ideal
class is not principal admits \genfar neighbours by Theorem
\ref{theo:bianverymax}, and they are not Farey neighbours since this
ideal is not principal. Therefore by Remark \ref{rem:remdefgenfar}
(\hyperlink{remdefgenfar1}{1}), there exist pairs of \genfar
neighbours that are not pairs of Farey neighbours if and only if
$h_K\geq 2$.

\medskip
\smallskip\noindent{\bf Proof of Theorem \ref{theo:bianverymax}.}
Since $\Ga_K$ preserves the set of pairs of \genfar neighbours and
since the stabilizer of any element of $\PP^1(K)$ is an infinite
parabolic subgroup, we only have to prove that every element
$x'\in\PP^1(K)$ admits an element $x$ in its $\Ga_K$-orbit that has a
\genfar neighbour $y$.

Let $a',b'\in\OOO_K$ be such that $x'=\frac{a'}{b'}$, and let $\aaa'=
a'\OOO_K +b'\OOO_K$. For all $a,b\in\OOO_K$, if $\aaa= a\OOO_K
+b\OOO_K$ belongs to the same ideal class as $\aaa'$, then $x=
\frac{a}{b}$ belongs to the same $\Ga_K$-orbit as $x'$ by the
bijection \eqref{cuspclass}. If $\aaa$ is principal, then $x$ is the
same $\Ga_K$-orbit as $\infty$, hence has \genfar neighbours by Remark
\ref{rem:remdefgenfar} (\hyperlink{remdefgenfar1}{1}), and so does
$x'$. The norm of a nonprincipal prime ideal is a prime integer. By
Weber's theorem in \cite[Sect.~X.12]{Cohn80}, there are infinitely
many prime ideals in each ideal class. By for instance
\cite[Thm.~6.14]{Lemmermeyer21}, we may hence assume that $\aaa$ is a
nonprincipal prime ideal such that $[\aaa]=[\aaa']$ and $\Nr(\aaa)=
p_0$ is an odd prime such that one of the following two claims
hold.\footnote{The third case of \cite[Thm.~6.14]{Lemmermeyer21} does
not occur, since otherwise $p_0\OOO_K$ would be a prime ideal in that
case, and $\Nr(\aaa)=p_0$ implies that $\aaa\mid p_0\OOO_K$, so that
$\aaa= p_0\OOO_K$ and $\Nr(\aaa)=p_0^{\,2}$, a contradiction.}

\medskip
\noindent {\bf Case \hypertarget{dichotomyi}{i)}. } The prime $p_0$
ramifies in $K$, that is $p_0\mid\disc_K$. With $\aaa=
\sqrt{\radix}\OOO_K +p_0 \OOO_K$, we have $p_0\OOO_K= \aaa^2$ by
loc. cit.. We define $a_0=0$ in Case \hyperlink{dichotomyi}{i}). Note
that we have $p_0\mid-\radix =\Nr(a_0+\sqrt{\radix})$ since $p_0$ is
odd (and $D_K$ and $f_K$ have the same odd prime factors), and
$p_0^2\not\mid -\radix= \Nr(a_0+\sqrt{\radix})$ since $\radix$ is
square-free.

\medskip
\noindent {\bf Case \hypertarget{dichotomyii}{ii)}. }  The prime
$p_0$ splits in $K$, that is the discriminant $D_K$ is a quadratic
residue modulo $p_0$. Since $p_0$ is odd, there exists $a_0\in \ZZ\ssm
p_0\ZZ$ such that $a_0^{\;2}=\radix\!\!\mod p_0$. Let us define $\aaa=
(a_0+ \sqrt{\radix}\,) \OOO_K +p_0 \OOO_K$. We then have $p_0\OOO_K= \aaa\,
\ov\aaa$ by loc.~cit..  We have $p_0\mid a_0^{\;2}- \radix = \Nr(a_0+
\sqrt{\radix}\,)$. If $p_0^{\;2}\mid \Nr(a_0+ \sqrt{\radix}\,)$, then
$p_0^{\;2}$ does not divide $\Nr(a_0+p_0+ \sqrt{\radix}\,)=\Nr(a_0+
\sqrt{\radix}\,)+p_0^{\;2}+2a_0p_0$ since $p_0$ is odd and $a_0\neq
0\!\!\mod p_0$. Hence up to replacing $a_0$ by $a_0+p_0$, which does
not change $\aaa$ nor the fact that $p_0\mid \Nr(a_0+
\sqrt{\radix}\,)$, we have $p_0^2\not\mid \Nr(a_0+\sqrt{\radix})$.

\medskip
In both cases, $\frac{\phantom{\underline f}\!\!\! \Nr(a_0+
  \sqrt{\radix} \,)}{p_0}$ and $p_0$ are relatively prime integers. By
Bézout's identity for $\ZZ$, there exist $t,u\in\ZZ$ such that
 \begin{equation}\label{eq:bezout}
 \frac{\Nr(a_0+\sqrt{\radix}\,)}{p_0}\, u-p_0\, t=1\,.
 \end{equation}
Thus, setting $a=a_0+\sqrt{\radix}$, $b=p_0$, $c=t(a_0+ \sqrt{\radix}
\,)$ and $d=u\,\frac{\phantom{\underline f}\!\!\!\Nr(a_0+\sqrt{\radix}
  \,)}{p_0}$ that all belong to $\OOO_K$, we have $\aaa= a\,\OOO_K+
b\,\OOO_K$ and $(a,b,c,d)$ satisfies Equation \eqref{eq:fareyideal} :
Using Equation \eqref{eq:bezout} for the last two equalities, we have
\begin{align}
&(a\OOO_K+b\OOO_K)(c\OOO_K+d\OOO_K)\nonumber
\\=\;& \big((a_0+\sqrt{\radix}\,)\OOO_K+p_0\,\OOO_K\big)
\big(t\,(a_0+\sqrt{\radix}\,)\OOO_K+
u\, \frac{\Nr(a_0+\sqrt{\radix}\,)}{p_0}\OOO_K\big)\nonumber\\
=\;&(a_0+\sqrt{\radix}\,)\Big(t\,(a_0+\sqrt{\radix}\,)\OOO_K
+u\,\frac{\Nr(a_0+\sqrt{\radix}\,)}{p_0}\OOO_K+p_0\,t\,\OOO_K
+u\, (a_0-\sqrt{\radix}\,)\OOO_K\Big)\nonumber\\=\;&
(a_0+\sqrt{\radix}\,)\OOO_K=(ad-bc)\OOO_K\,.\label{eq:computdemtheo0}
\end{align}
Thus $x=\frac ab$ is in the same $\Ga_K$-orbit as $x'$ since they have
the same associated ideal class $[\aaa]=[\aaa']$, and $y=\frac cd$ is
a \genfar neighbour of $x$ as wanted.
\cqfd

\medskip
Note that the computation \eqref{eq:computdemtheo0} in the proof of
Theorem \ref{theo:bianverymax} is valid as long as the integers
$\frac{\phantom{\underline f}\!\!\!\Nr(a_0+\sqrt{\radix}\,)} {p_0}$
and $p_0$ in Equation \eqref{eq:bezout} are relatively prime.  Thus,
in order to produce examples of \genfar neighbours, we can use the
tables at the end of \cite{Sommer11} (and their reproduction at the
end of \cite{Cohn80}) where representatives are listed for all ideal
classes of imaginary quadratic number fields with $-97\le \radix\le
-1$. Note that the representatives in these tables are not always
prime ideals, though.

Counting results for pairs of \genfar neighbours in an orbit of a
given pair $\{x,y\}$ by a finite index subgroup $\Ga$ of $\Ga_K$
follow immediately from Theorem \ref{theo:iqcount}. The results become
more explicit in the cases where the values of the reciprocity index
$\iota_\Ga(x,y)$ and the multiplicity $m_\Ga(x,y)$ are known.  If the
\genfar neighbours $x,y\in K$ are in two different $\Ga_K$-orbits,
then $\iota_{\Ga_K}(x,y)=2$. The following examples provide in
particular infinite collections of \genfar neighbours $x,y\in K$ with
$\iota_{\Ga_K}(x,y)= 1$.

\bexems\label{ex:rami}
\hypertarget{exampfarid1}{(1)} Assume that $-\radix$ is at least $6$
and is not a prime. Let $p_0$ be a prime factor of $-\radix$. Then
$p_0$ is ramified in $K$. As in Case \hyperlink{dichotomyi} {i}) in
the proof of Theorem \ref{theo:bianverymax}, the ideal $\aaa=
\sqrt{\radix}\OOO_K+ p_0\OOO_K$ satisfies $\aaa^2=p_0\OOO_K$. In
particular $\Nr(\aaa)=p_0$, hence $\aaa$ is prime.  Furthermore $\aaa$
is not principal by the following result.

\blemm
There are no elements of norm $p_0$ in $\OOO_K$.
\elemm

\dem Let $\omega_K=\sqrt{\radix}$ if $\radix\equiv 2,3\bmod 4$ and
$\omega_K=\frac{\phantom{\underline f}\!\!\!1+\sqrt{\radix}}2$ if
$\radix\equiv 1\bmod 4$, so that $\OOO_K=\ZZ+\ZZ\, \omega_K$.  The
norm of an element of $\OOO_K\cap\ZZ$ is not a prime. Let $u,v \in\ZZ$
with $v\neq 0$, and assume for a contradiction that
$\Nr(u+v\,\omega_K)=p_0$. If $\radix\not\equiv 1\bmod 4$, then
$\Nr(u+v\,\omega_K)\geq -v^2\radix \geq-\radix>p_0$ since $p_0\mid
-\radix$ and $\radix$ is not a prime, a contradiction.  If $|v|\geq 2$
or if $\frac{\phantom{\underline f} \!\!\!-\radix}{p_0}>4$, then
$\Nr(u+v\,\omega_K)\geq v^2(\Im\; \omega_K)^2 \geq -v^2\frac{\radix}4
>p_0$, a contradiction. Thus $\radix\equiv 1\bmod 4$ and in particular
$\radix$ is odd with $-\radix =3\,p_0$, and $v= \pm1$. Therefore
$\Nr(u+v\,\omega_K)= (u\pm\frac{1}{2})^2 +\frac{3\,p_0}{4}$. Since the
solutions of the equations $(u\pm \frac{1}{2})^2 +\frac{3\,p_0}{4}=p_0$
with unknown $u$, that are $u=\frac{\pm1\pm\sqrt{p_0}}{2}$, are
irrational, this contradicts the assumption.
\cqfd

\medskip
Equation \eqref{eq:bezout} becomes
\begin{equation}\label{eq:bezoutex}
-\frac{\radix}{p_0} \,u-p_0\,t=1\,.
\end{equation}
Let $(t,u)\in\ZZ^2$ be an integral solution of Equation
\eqref{eq:bezoutex}.  By the end of the proof of Theorem
\ref{theo:bianverymax}, the element $\beta = \frac{\phantom{
\underline f}\!\!\!t\,\sqrt{\radix}}{- \phantom{\hat{\hat f}}\!\!\! u
\frac{\radix}{p_0}}=-\frac{\phantom{\underline f}\!\!\!t\; p_0}
{\phantom{\hat f} \!\!\!u \sqrt{\radix}} \in K$ is a \genfar neighbour
of $\alpha= \frac{\phantom{\underline f} \!\!\!  \sqrt{\radix}} {p_0}$
(but note that they are not Farey neighbours since $\aaa$ is not
principal).  By Equation \eqref{eq:bezoutex}, we have
$C=\begin{bmatrix}\sqrt{\radix}&t\\p_0&-\frac{\phantom{\underline f}
  \!\!\!u\sqrt{\radix}}{p_0} \end{bmatrix}\in \PSL_2(K)$. We have
$C\cdot \infty=\alpha$ and $C\cdot 0=\beta$. Note that $\iota_{p_0}
=\begin{bsmallmatrix} 0 &\frac{1}{p_0} \\-p_0 &0\end{bsmallmatrix}
\in\PSL_2(K)$ is an involution exchanging the points at infinity
$\infty$ and $0$ of $\partial_\infty\htr$. The conjugate involution
\[
E=C\;\iota_{p_0}\;C^{-1}=\begin{bmatrix}(tu-1)\sqrt{\radix} &
t^2\,p_0+\frac{\phantom{ \underline f} \!\!\!\radix}{p_0}\\[3pt]
-\frac{\phantom{\underline f} \!\!\!u^2\,\radix}{p_0}-p_0&
(1-tu)\sqrt{\radix} \end{bmatrix}
\]
belongs to $\PSLOK$ since $p_0\mid -f_K$ and satisfies $E\cdot\alpha=
(C\;\iota_{p_0})\cdot \infty=C\cdot 0=\beta$ and similarly $E\cdot
\beta=\alpha$. Thus, with the notation of Section \ref{sec:commperp},
the divergent geodesic $\ell_{\alpha,\beta}$ in $M_K$ is reciprocal
and $\iota_{\Ga_K,\,{\rm rec}} (\alpha,\beta)=1$.

The pointwise stabilizer $\Ga_{\alpha,\beta}$ of $\wt\ell_{\alpha,
\beta}(\RR)=C\cdot\wt\ell_{\infty, 0}(\RR)$, which is the conjugate by
$C$ of the pointwise stabilizer $\Big\{M(\theta)= \begin{bsmallmatrix}
  e^{i\theta}&0 \\ 0&e^{-i\theta} \end{bsmallmatrix} :\theta\in\RR
\Big\}$ of $\wt\ell_{\infty, 0}(\RR)$, can also be determined. Note
that by Equation \eqref{eq:bezoutex}, we have
\begin{align*}
C\,M(\theta)\,C^{-1} &
=\begin{bmatrix} \cos\theta-i\,(1+\frac{2\radix u}{p_0})\sin \theta &
-2i\,t\,\sqrt{\radix}\sin\theta\\ -2i\,u\,\sqrt{\radix}\sin\theta&
\cos\theta+i\,(1+\frac{2\radix u}{p_0})\sin \theta
\end{bmatrix}\,.
\end{align*}
Let $\theta\in\RR$ be such that $C\,M(\theta)\,C^{-1}$ belongs to
$\PSLOK$. Then the trace of this matrix, which is $\pm2\cos \theta$,
belongs to $\OOO_K\cap\RR=\ZZ$. Therefore $\cos\theta= 0,\pm 1,
\pm\frac 12$ and correspondingly $\sin\theta=\pm 1,0,\pm
\frac{\sqrt{3}}{2}$. If, for a contradiction, $\sin\theta\neq 0$, then
the $2$-$1$ entry of the above matrix, that is equal to $\pm\, 2\,u\,
\sqrt{-\radix}$ or $\pm\, u\,\sqrt{-3\,\radix}$, also belongs to
$\OOO_K\cap\RR =\ZZ$. But since $-\radix$ is squarefree and at least
$6>3$, these entries are irrational, a contradiction. Thus, the
stabilizer $\Ga_{\alpha,\beta}$ is trivial and $m_{\Ga_K}(\alpha,
\beta) =1$.

\medskip\noindent
\hypertarget{exampfarid2}{(2)} Assume in this family of examples that
$K=\QQ(\sqrt{\radix}\,)$ with $\radix\equiv 3\bmod 4$ and $-\radix\geq
5$. Then $\OOO_K=\ZZ+\sqrt{\radix}\,\ZZ$, $D_K=4\, \radix$ and $2\mid
D_K$, so that $2$ ramifies in $K$. Let $\aaa=(1+\sqrt{\radix}\,)\OOO_K
+2\OOO_K$. Since $\overline{1+\sqrt{\radix}\,} =2-(1+ \sqrt{\radix}
\,)$ and by for instance \cite[Lem.~13.8.4]{Artin91}, we have $\aaa^2=
\aaa\,\overline{\aaa}= 2\OOO_K$, so that the class of $\aaa$ has order
$2$ in the ideal class group $\icg K$.  Even if the case $p_0=2$ does
not appear in the proof of Theorem \ref{theo:bianverymax}, the
analogous computations work in the present infinite collection of
examples.

Now, $\Nr(1+\sqrt{\radix}\,)=1-\radix\equiv 2\bmod 4$. Hence
$\frac{\phantom{\underline f}\!\!\!1-\radix}{2}$ and $2$ are coprime
integers, and Equation \eqref{eq:bezout} becomes
$\frac{\phantom{\underline f}\!\!\!  1-\radix}2\, u-2\,t=1$, satisfied
for instance by $t=- \frac{\phantom{\underline f}\!\!\!1+\radix}4$ and
$u=1$. Let $a=1+\sqrt{\radix}$, $b=2$, $c=-\frac{\phantom{\underline
    f} \!\!\!1+ \radix} {2}$ and $d=1-\sqrt{\radix}$, that belong to
$\OOO_K$ and satisfy $ad-bc=2$. By a computation similar to the one in
Equation \eqref{eq:computdemtheo0}, the quadruple $(a,b,c,d)$ is a
solution of Equation \eqref{eq:fareyideal}. Hence $\alpha=\frac
ab=\frac{\phantom{\underline f}\!\!\!1+\sqrt{\radix}} {2}$ and
$\beta=\frac cd= -\frac{\phantom{\underline f}\!\!\!  1+\radix}
{\phantom{\hat f}\!\!\!2(1- \sqrt{\radix}\,)}$ are \genfar neighbours
(that are not Farey neighbours since $\aaa$ is not principal).

The element $C=\begin{bmatrix} 1+\sqrt{\radix}
&-\frac{\phantom{\underline f}\!\!\!1+\radix}4 \\[3pt]2
&\frac{\phantom{\underline f}\!\!\!1- \sqrt{\radix}}{2} \end{bmatrix}
\in\PSL_2(K)$ maps $\infty$ and $0$ to $\alpha$ and $\beta$
respectively. Note that $\iota_{2}=\begin{bsmallmatrix} \;0
&\frac{1}{2} \\-2 &0\end{bsmallmatrix}\in\PSL_2(K)$ is an involution
exchanging $\infty$ and $0$. The involution
\[
E=C\;\iota_2\;C^{-1}=\begin{bmatrix}
\frac{\phantom{\underline f}\!\!\!-3+\radix-(5+\radix)\sqrt{\radix}}4 &
\frac{\phantom{\underline f}\!\!\!5+6\radix+\radix^2}8+\sqrt{\radix}
\\[3pt]-\frac{\phantom{\underline f}\!\!\!5+\radix}2+\sqrt{\radix}&
 \frac{\phantom{\underline f}\!\!\!3-\radix+(5+\radix)\sqrt{\radix}}4
\end{bmatrix}
\]
belongs to $\PSLOK$ since $\radix\equiv 3\bmod 4$, and satisfies
$E\cdot\alpha=\beta$ and $E\cdot \beta=\alpha$. Thus, the divergent
geodesic $\ell_{\alpha,\beta}$ in $M_K$ is reciprocal and
$\iota_{\Ga_K,\,{\rm rec}} (\alpha, \beta) =1$.

The pointwise stabilizer $\Ga_{\alpha,\beta}$ of the geodesic line
$\wt\ell_{\alpha,\beta}(\RR)= C\cdot\wt\ell_{\infty,0}(\RR)$ can be
determined as in the previous examples (\hyperlink{exampfarid1}{1}).
Let $\theta\in\RR$ be such that the entries of the elliptic element
\[
C\;M(\theta)\;C^{-1}=
\begin{bmatrix} \cos\theta-i\radix\sin\theta&
  i(1+\sqrt{\radix})\frac{\phantom{\underline f}\!\!\!1+\radix}2\sin\theta
  \\[3pt] -2i(\sqrt{\radix}-1)\sin\theta&\cos\theta+i\radix\sin\theta
\end{bmatrix}
\]
are in $\OOO_K$. Then by taking the sum and the difference of the
diagonal entries, we have $2\cos\theta\in\OOO_K\cap \RR=\ZZ$ and $2i
\radix\sin\theta\in\OOO_K\cap(i\RR)= \sqrt{\radix}\,\ZZ$. Hence $\cos
\theta =0,\pm1,\pm\frac 12$ and $2\,\sqrt{-\radix}\sin\theta\in\ZZ$.
Since $-\radix\ge 5>3 $, this implies as in the previous examples
(\hyperlink{exampfarid1}{1}) that the stabilizer $\Ga_{\alpha,\beta}$
is trivial, so that $m_{\Ga_K}(\alpha,\beta) =1$.  Note that if
$K=\QQ(\sqrt{-1})$, then $\Ga_{\alpha,\beta}$ consists of $\id$ and
$\begin{bmatrix} i&0\\2(1+i)&-i\end{bmatrix}$, so that
  $m_{\Ga_K}(\alpha,\beta) =2$.

\medskip
\noindent\hypertarget{exampfarid2b}{(2)$^{\rm bis}$} Let us consider
the particular case $\radix=-5$ of the previous family of examples
(\hyperlink{exampfarid2}{2}). Recall that the class number of $K=
\QQ(\sqrt{- 5})$ is $h_K=2$. The ideal $\aaa=(1+i\sqrt 5)\OOO_K
+2\OOO_K$ (that satisfies $\aaa=\overline{\aaa}$ and $\Nr(\aaa)=2$) is
a prime representative of the unique nonprincipal ideal class. By the
general computation above, the elements $\alpha=\frac{1+i\sqrt 5}2$
and $\beta=-\frac 2{1-i\sqrt 5}$ in $\QQ(\sqrt{- 5})$ are \genfar
neighbours, that are not Farey neighbours.

The orbit of $\beta$ under the stabilizer $\Ga_{\alpha}$ of $\alpha$
in $\Ga_K$ gives an infinite collection of \genfar  neighbours of
$\alpha$. For example by the proof of \cite[Lemma 6]{ParPau11BLMS},
we have
\[
\Ga_{\alpha}=\bigg\{\begin{bmatrix} 1+(1+i\sqrt 5) x & (2-i\sqrt 5)x\\
2x& 1-(1+i\sqrt 5) x\end{bmatrix} : x\in\OOO_K\bigg\}\,.
\]
The figures below show the canonical horoball $B_\alpha$ (drawn in
red) and the canonical horoballs (drawn in beige) of some of the
\genfar neighbours of $\alpha$ (images of $B_\beta$ by elements of
$\Ga_\alpha$), that are hence tangent.

\begin{center}
\includegraphics{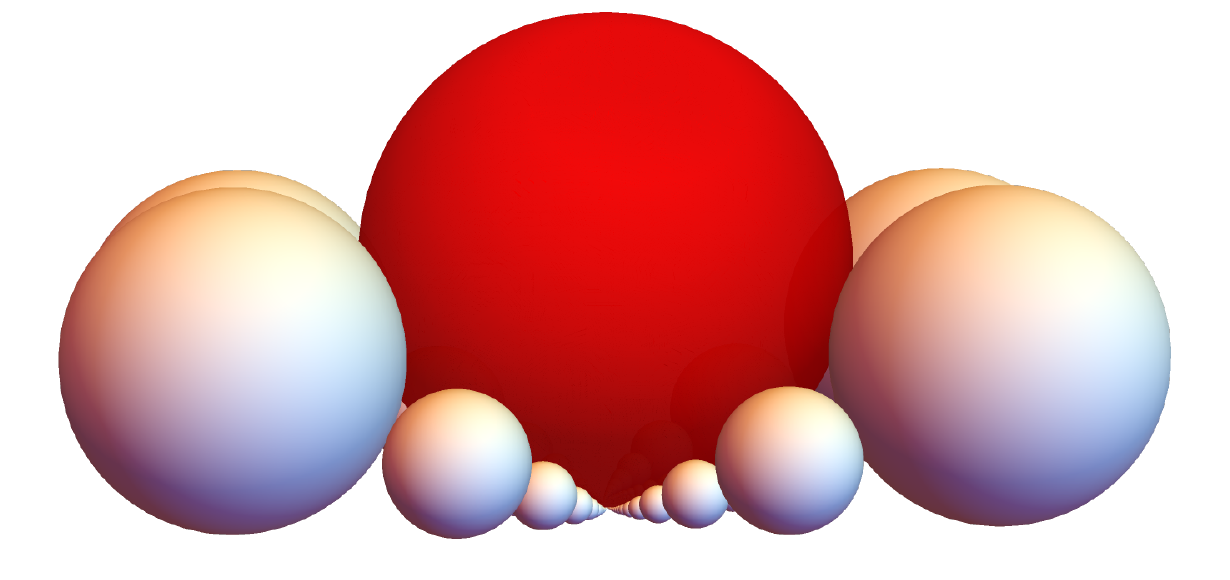}
\includegraphics[width=.8\textwidth]{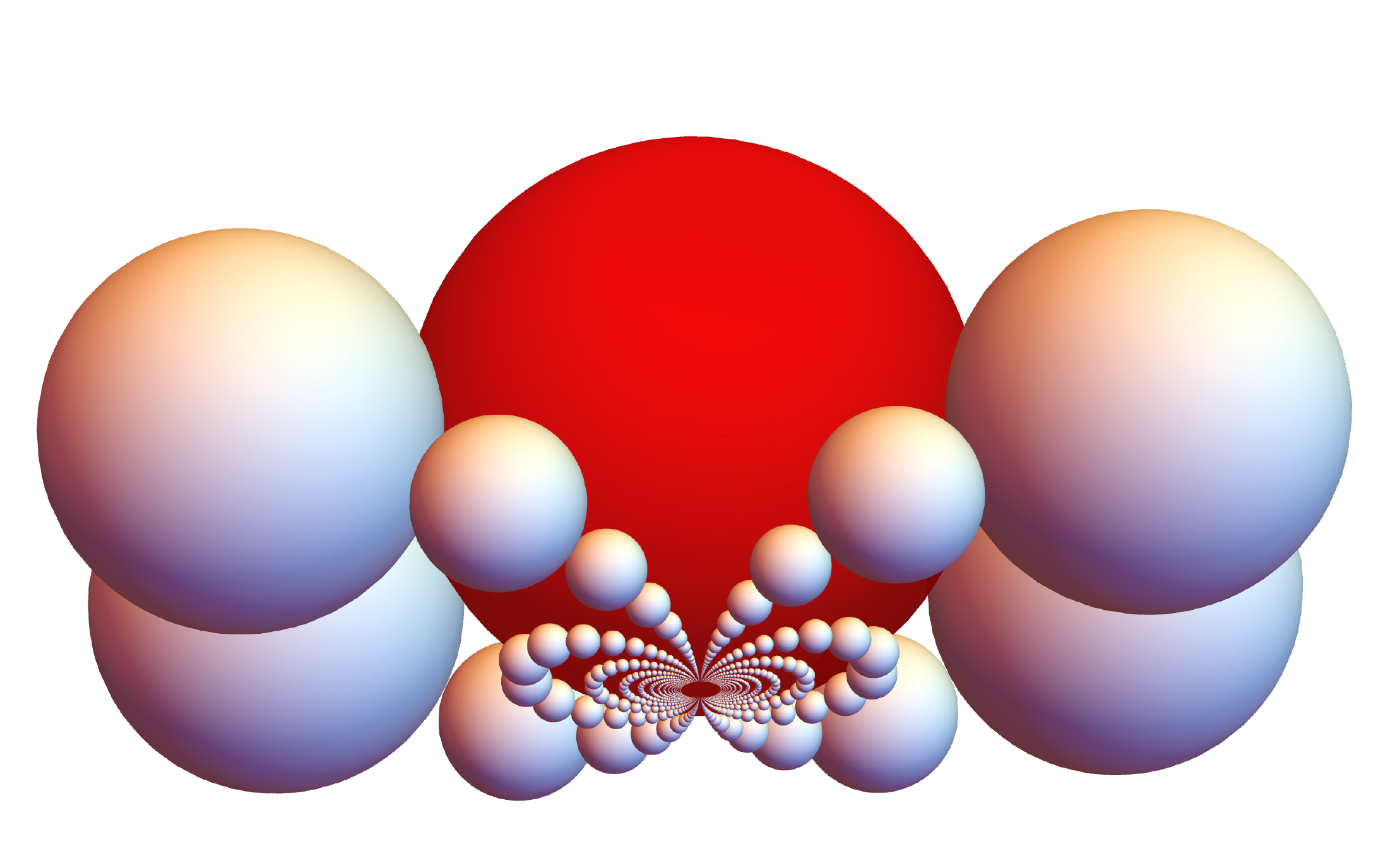}
\end{center}

As seen above, the elliptic element $\begin{bmatrix} -2&i\sqrt{5}
\\ \ i\sqrt{5}& 2\end{bmatrix}\in\PSLOK$ of order $2$ exchanges
$\alpha$ and $\beta$. Hence for all \genfar  neighbours $\alpha'$ and
$\beta'$ in $K$, the divergent geodesic $\ell_{\alpha',\beta'}$ in
$M_K$ is reciprocal.

The figure below shows parts of the two families of horospheres that
correspond to the classes of principal (blue) and non-principal
(orange) ideals.  The horospheres are somewhat translucent, and the
horospheres can be seen even if they are behind other canonical
horoballs as seen from the viewpoint.  The horospheres are restricted
to the symmetric closed fundamental domain $\{(z,t)\in\htr: |\Re z|\le
\frac 12,\ |\Im z|\le\frac{\sqrt 5}2\}$ of the stabilizer of $\infty$
in $\PSLOK$, and the picture is cut at height $\frac 13$.

\begin{center}
\includegraphics[width=.8\textwidth]{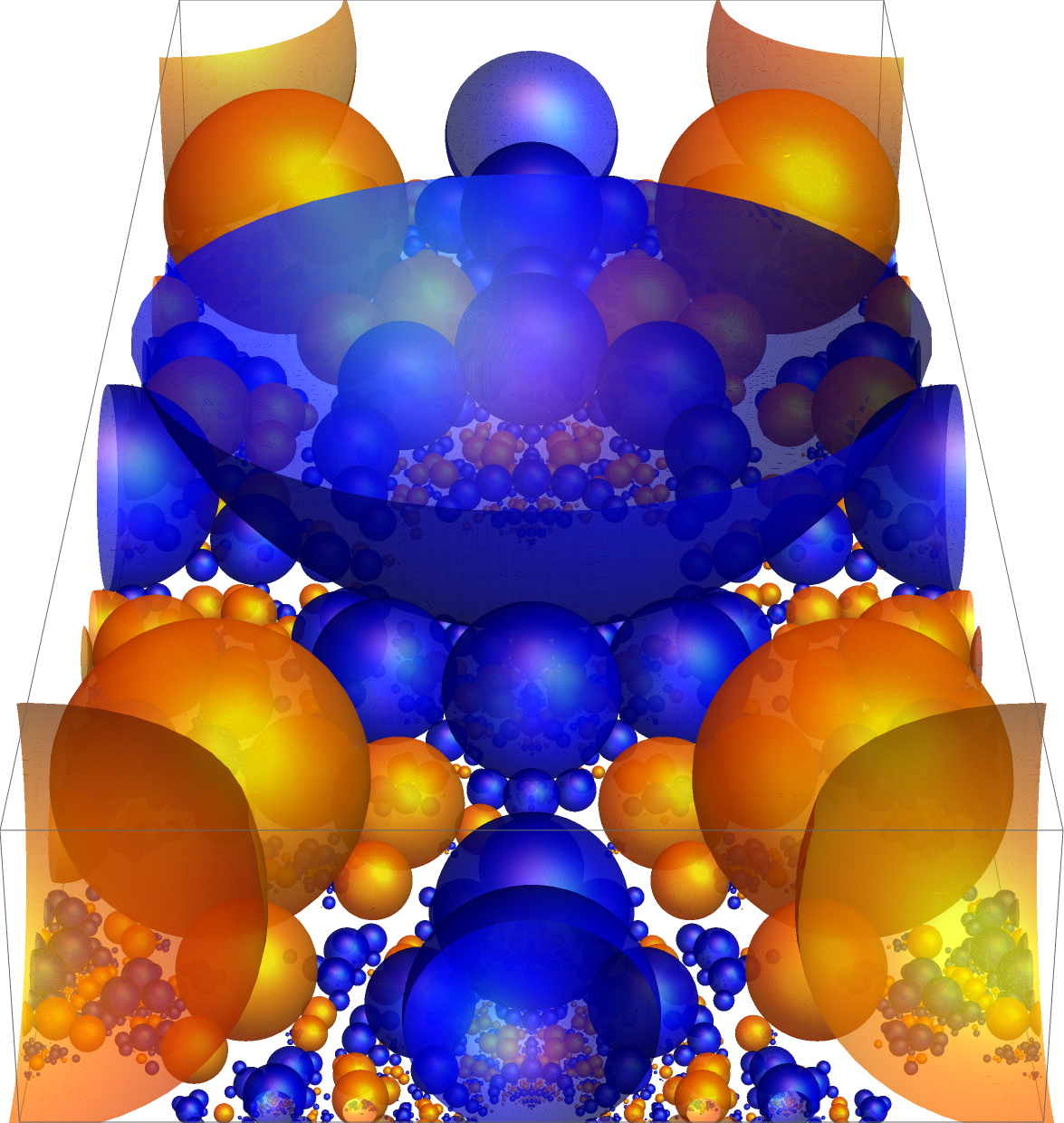}
\end{center}

\medskip
\noindent\hypertarget{exampfarid3}{(3)} Assume in this other family of
examples that $K=\QQ(\sqrt{\radix}\,)$ with $\radix\equiv 2\bmod 4$
and $-\radix\ge 6$. Similarly as for the examples
(\hyperlink{exampfarid2}{2}), let $a= \sqrt{\radix}$, $b=2$, $c=
\frac{\phantom{\underline f}\!\!\!\radix+2}2$ and $d=\sqrt{\radix}$,
that belong to $\OOO_K$.  The nonprincipal ideal class of
$\aaa=a\OOO_K+b\OOO_K$ has order $2$ in $\icg K$ again by for instance
\cite[Lem.~13.8.4]{Artin91} since $2$ ramifies in $K$. Then we have
$ad-bc= -2$ and $(a,b,c,d)$ satisfies the condition
\eqref{eq:fareyideal}. Thus the elements $\alpha= \frac ab
=\frac{\phantom{\underline f}\!\!\!\sqrt{\radix}}2$ and
$\beta=\frac{\phantom{\underline f}\!\!\!2+ \radix}
{\phantom{\underline f}\!\!\!2\sqrt{\radix}}$ of $K$ are
\genfar  neighbours that are not Farey neighbours (since $[\aaa]$ has
order $2$ hence is not the principal class). The element
$C=\begin{bmatrix} \sqrt{\radix} &-\frac{\phantom{\underline f}\!\!\!2
  +\radix}4\\[3pt] 2 &-\frac{\phantom{\underline f}\!\!\!
  \sqrt{\radix}}{2}\end{bmatrix}\in\PSL_2(K)$ maps $\infty$ and $0$ to
$\alpha$ and $\beta$ respectively. The involution
\[
E=C\;\iota_2\;C^{-1}=
\begin{bmatrix} 
-\frac{\phantom{\underline f}\!\!\!6+\radix}4\sqrt{\radix} &
\frac{\phantom{\underline f}\!\!\!4+8\radix+\radix^2}8\\[3pt]
-\frac{\phantom{\underline f}\!\!\!4+\radix}2& \frac{\phantom{\underline
    f}\!\!\!6+\radix}4\sqrt{\radix}
\end{bmatrix}
\]
belongs to $\PSLOK$ and satisfies $E\cdot\alpha=\beta$ and $E\cdot
\beta =\alpha$. Thus, the divergent geodesic $\ell_{\alpha,\beta}$ in
$M_K$ is reciprocal and $\iota_{\Ga_K,\,{\rm rec}} (\alpha,\beta) =1$.

The pointwise stabilizer $\Ga_{\alpha,\beta}$ of the geodesic line
$\wt\ell_{\alpha,\beta}(\RR)$ can be determined in the same way. For
$\theta\in\RR$, the entries of the elliptic element
\[
C\;M(\theta)\; C^{-1}=
\begin{bmatrix} 
  \cos\theta-i(1+\radix)\sin\theta&
  i\sqrt{\radix}\,\frac{\phantom{\underline f}2+\radix}2\sin\theta
\\[3pt]-2i\sqrt{\radix}\sin\theta&\cos\theta+i(1+\radix)\sin\theta
\end{bmatrix}
\]
are in $\OOO_K$ if and only if $\sin\theta=0$ (as for the examples
(\hyperlink{exampfarid2}{2}), otherwise the $2$-$1$ entry
would be an irrational real number), since $-\radix \geq 6>3$.  This
implies that $\Ga_{\alpha,\beta}$ is trivial, thus
$m_{\Ga_K}(\alpha,\beta) =1$.

\medskip
\noindent\hypertarget{exampfarid4}{(4)} The last congruence property
on $\radix$ is when $\radix\equiv 1\bmod 4$. If furthermore $-\radix$
is a prime integer, then there are no elements of order $2$ in the
class group $\icg K$ by, for instance, \cite[Prop.~3.11]{Cox13}. By
Remark \ref{rem:remdefgenfar} (\hyperlink{remdefgenfar2}{2}), all
reciprocal \genfar  neighbours in $K$ are then Farey neighbours.
\eexems

\section{Farey neighbours in rational definite quaternion
  algebras}\label{sec:quaternion}

In this Section, we study similar asymptotic countings of quaternionic
Farey neighbours. Let $\HH$ be the standard Hamilton quaternion
algebra over $\RR$, with canonical $\RR$-basis $(1,i,j,k)$ and with
conjugation $x\mapsto \overline{x}$, reduced norm $\n$ and reduced
trace $\tr$. We denote by $\PP^1_r(\HH)$ the right projective line of
$\HH$, identified as usual with the Alexandrov compactification
$\HH\cup\{\infty\}$ where $[x:y]=xy^{-1}$ if $y\neq 0$ and $[1:0] =1
\;0^{-1}=\infty$. 

Let $\OOO$ be a maximal order in a quaternion
algebra $A$ over $\QQ$, which is definite (that is, $A\otimes_\QQ
\RR=\HH$), with class number $h_A$ and discriminant $D_A$. Its group
$\OOO^\times$ of invertible elements is finite, of order $2$, $4$,
$6$, $12$ (when $D_A=3$) or $24$ (when $D_A=2$). An example is given
by the {\it Hurwitz order} $\OOO=\ZZ+\ZZ\, i+\ZZ\, j+
\ZZ\,\frac{1+i+j+k}{2}$ in $A=\QQ+\QQ\, i + \QQ\, j + \QQ\, k$, in
which case $h_A=1$ and $D_A=2$. We refer for these informations and
more to \cite{Vigneras80}.

We will say that two elements $\alpha$ and $\beta$ in $\PP^1_r(A)=
A\cup\{\infty\}$ are {\it Farey neighbours with respect to}
$\OOO$ if there exist $p,q,r,s\in \OOO$ with $\alpha=pq^{-1}$,
$\beta=rs^{-1}$, and either we have $q=0$ and $p,s\in \OOO^\times$ or
we have $q\neq 0$ and
\begin{equation}\label{eq:defiQuatFarNei}
\n(qpq^{-1}s-qr)=1\,.
\end{equation}
This condition is the appropriate noncommutative analog of Equation
\eqref{eq:defiFareyNeigh}. Let $\fn_\OOO$ be the set of unordered
pairs of Farey neighbours in $\PP^1_r(A)$ with respect to $\OOO$. It
is easy to check that the additive group $\OOO$ acts by simultaneous
translations on the set $\fn_\OOO$. The following theorem gives an
effective asymptotic counting result for pairs of quaternionic Farey
neighbours with respect to $\OOO$ when the lower bound on their
distances shrinks to $0$.

\btheo\label{theo:quaternionfareycount} As $\epsilon>0$ tends to $0$,
we have
\begin{align*}
&\card\big(\OOO\bs\big\{ \{\alpha,\beta\}\in \fn_\OOO :
\n(\beta-\alpha) \ge \epsilon\big\}\big) \\=\;&-\frac{2160\;D_A}
{\zeta(3)\;|\OOO^\times|^2\;\prod_{p\,\mid D_A}(p^3-1)(p-1)}\;
\frac{\ln\epsilon}{\epsilon^2}+\bigO(\epsilon^{-2})\, .
\end{align*}
\etheo

As usual, the above index $p$ ranges over primes. We will actually
prove a much stronger result, that requires some information on the
Hamil\-ton-Bian\-chi groups $\PSL_2(\OOO)$. See for instance
\cite{Kellerhals03} for background; we will follow the presentation of
\cite[\S 3]{ParPau13ANT}.

The {\it Dieudonné determinant} is the group morphism $\Det :
\GL_2(\HH)\ra\;]0,+\infty[$ defined by
\begin{equation}\label{eq:detdieu}
 \Det:\begin{pmatrix} a & b \\ c & d\end{pmatrix} \mapsto
   \big(\n(a\, d) + \n(b\, c)- \tr(a\, \ov{c}\, d\,
   \ov{b}\,)\big)^{\frac 12}\,.
\end{equation}
The Lie group $\SL_2(\HH)$ is the kernel of $\Det$. We denote by
$\begin{bsmallmatrix} a & b \\ c & d\end{bsmallmatrix}\in \PSL_2(\HH)
=\SL_2(\HH)/\{\pm\id\}$ the image of $\begin{psmallmatrix}a & b \\ c
& d\end{psmallmatrix}\in\SL_2(\HH)$.  The group $\PSL_2(\HH)$ acts
faithfully by homographies on the right projective plane
$\PP^1_r(\HH)=\HH\cup \{\infty\}$, by $\begin{bsmallmatrix} a & b \\ c
  &d\end{bsmallmatrix} \cdot z=(az+b)(cz+d)^{-1}$ for all
$\begin{bsmallmatrix}a & b \\ c & d\end{bsmallmatrix}\in\PSL_2(\HH)$
and $z\in\PP^1_r(\HH)$, with the usual conventions when $z=\infty$ or
$z=-c^{-1}d$. With $ds^2_\HH$ the usual translation-invariant
flat Riemannian metric on $\HH$ (making the canonical $\RR$-basis
$(1,i,j,k)$ of $\HH$ orthonormal at each point), we identify
$\HH^5_\RR$ with
\[
\Big(\{(z,t)\in\HH\times\RR:t>0\},\;
\frac{ds^2_\HH+dt^2}{t^2} \Big)\,.
\]
The group $\PSL_2(\HH)$ acts faithfully on $\HH^5_\RR$ by the Poincaré
extension procedure,\footnote{See for instance
\cite[Eq.~(14)]{ParPau13ANT}} and $\PSL_2(\HH)$ thus identifies with
the orientation preserving isometry group of $\HH^5_\RR$.

Let $_\OOO\I$ be the set of left ideal classes of $\OOO$, whose cardinality 
is the class number $h_A$. The subgroup $\Ga_\OOO= \PSL_2(\OOO)$ is an
arithmetic lattice in $\PSL_2(\HH)$.  By for instance \cite[Satz 2.1,
  2.2] {KraOse90}, it acts with $(h_A)^{2}$ orbits on its set of
parabolic fixed points $\PP^1_r(A)=A\cup \{\infty\}$ in
$\partial_\infty\hcr$. The stabiliser of $\infty$ in $\Ga_\OOO$ is
\[
\Ga_{\OOO,\infty}=\bigg\{\begin{bmatrix} a & b\\0&d\end{bmatrix}\in
\Ga_\OOO: a,d\in\OOO^\times\,,\; b\in\OOO\bigg\}\,.
\]

\btheo\label{theo:quatercount} Let $\Ga$ be a finite index subgroup of
$\Ga_\OOO=\PSL_2(\OOO)$, and let $\Ga_\infty$ be the stabiliser of
$\infty$ in $\Ga$.  For all distinct $x,y\in A\cup\{\infty\}$, as
$\epsilon>0$ tends to $0$, we have
\begin{align*}
&\card\big(\Ga_{\infty}\bs\big\{ \{\alpha,\beta\}\in \Ga\cdot\{x,y\} :
\n(\beta-\alpha)\ge\epsilon\big\}\big)\\=&
-\frac{2160\;D_A\;\iota_{\Ga,\,{\rm rec}}(x,y)\;
  [\Ga_{\OOO,\infty}:\Ga_\infty]}{\zeta(3)\;|\OOO^\times|^2\;
  \big(\prod_{p\,\mid D_A}(p^3-1)(p-1) \big)\;m_\Ga(x,y)\;[\Ga_{\OOO}:\Ga]}\;
\frac{\ln \epsilon}{\epsilon^2}+\bigO(\epsilon^{-2}).
\end{align*}
\etheo

\dem As in the proof of Theorem \ref{theo:rationalcount}, we apply
Theorem \ref{theo:1geoddivreal} with $n=5$, with $M=\Ga\bs\HH^5_\RR$,
with $D^-=\Ga_{\infty}\bs B_\infty$, and with $D^+=\ell_{x,y}(\RR)$.
By Emery's volume formula \cite[Theo.~8 and Appendix]{ParPau13ANT}, we
have
\[
\Vol(M)=[\Ga_\OOO:\Ga]\,\Vol(\Ga_\OOO\bs\HH^5_\RR)
=[\Ga_\OOO:\Ga]\;\frac{\zeta(3)\;\prod_{p\,\mid D_A}(p^3-1)(p-1)}
{11520}\,.
\]
The index $[\Ga_{\OOO,\infty}:\OOO]$ in $\Ga_{\OOO,\infty}$ of its
unipotent subgroup consisting in the translations by elements of
$\OOO$ is equal to $\frac{|\OOO^\times|^2}{2}$. By the Remark just
above \cite[Lemma 15]{ParPau13ANT}, we have
\[
\Vol(\partial D^-)=[\Ga_{\OOO,\infty}:\Ga_\infty]
\Vol(\Ga_{\OOO,\infty}\bs \partial B_\infty)=
[\Ga_{\OOO,\infty}:\Ga_\infty]\;\frac{D_A}{8\;|\OOO^\times|^2}\,.
\]
The Euclidean distance in $\HH$ between two elements
$\alpha,\beta\in\HH$ is $\n(\beta-\alpha)^{\frac{1}{2}}$, so that the
length of the common perpendicular from $B_\infty$ to
$\wt\ell_{\alpha,\beta} (\RR)$ when it exists is
$\ln\Big(\frac2{\n(\beta-\alpha)^{\frac{1}{2}}}\Big)$.  Hence as in
the proof of Theorem \ref{theo:rationalcount}, since
$\Ga\big(\frac{5}{2}\big)=\frac{3\sqrt{\pi}}{4}$ and $\Ga(3)=2!=2$, as
$\epsilon>0$ tends to $0$, we have
\begin{align*}
& \card\big(\Ga_\infty\bs\big\{ \{\alpha,\beta\}\in \Ga\cdot\{x,y\}:\
  \n(\beta-\alpha) \ge \epsilon\big\}\big)=
  \N_{D^-,D^+}\Big(\ln\frac 2{\sqrt{\epsilon}}\Big)+
  \bigO\Big(\ln\frac 2{\sqrt{\epsilon}}\Big)\\=&\,
  \frac{\Ga(\frac{5}{2})\;\iota_{\Ga,\,{\rm rec}}(x,y)\;
    [\Ga_{\OOO,\infty}:\Ga_\infty]\;\frac{D_A}
    {8\,|\OOO^\times|^2}}{2\,\sqrt{\pi}\;\Ga(3)\;m_\Ga(x,y)\;
    [\Ga_\OOO:\Ga]\,\frac{\zeta(3)\;\prod_{p\,\mid D_A}(p^3-1)(p-1)}
{11520}}\; \Big(\ln\frac 2{\sqrt{\epsilon}}\Big)
  \Big(\frac 2{\sqrt{\epsilon}}\Big)^4 + \bigO(\epsilon^{-2})\\= &
  -\frac{2160\;D_A\;\iota_{\Ga,\,{\rm rec}}(x,y)\;
  [\Ga_{\OOO,\infty}:\Ga_\infty]}{\zeta(3)\;|\OOO^\times|^2\;
  \big(\prod_{p\,\mid D_A}(p^3-1)(p-1) \big)\;m_\Ga(x,y)\;[\Ga_{\OOO}:\Ga]}\;
\frac{\ln \epsilon}{\epsilon^2}+\bigO(\epsilon^{-2})\,.\qquad\quad\Box
\end{align*}

In order to prove Theorem \ref{theo:quaternionfareycount}, the key
translation between the arithmetics and the geometry is the following
lemma.

\blemm\label{lem:fromarithtogeomquat} Two distinct elements $\alpha,
\beta\in\PP_r^1(A)=A\cup\{\infty\}$ are Farey neighbours with respect
to $\OOO$ if and only if there exists $\ga\in\Ga_\OOO= \PSL_2(\OOO)$
such that $\ga\cdot \infty= \alpha$ and $\ga\cdot 0= \beta$.
\elemm

\dem
For every $\begin{psmallmatrix} a & b \\ c & d\end{psmallmatrix} \in
\GL_2(\HH)$ such that $c\neq 0$, by for instance
\cite[Eq.~(12)]{ParPau13ANT}, we have $\Det\begin{psmallmatrix} a & b
\\ c & d\end{psmallmatrix}=\n(cac^{-1}d-cb)^{\frac{1}{2}}$.

Let $\alpha,\beta\in \PP^1_r(A)$ be distinct elements. Assume that
there exists $\ga=\begin{bsmallmatrix} p&r\\q&s\end{bsmallmatrix} \in
\Ga_\OOO$ such that $\ga\cdot \infty= \alpha$ and $\ga\cdot 0= \beta$.
If $q=0$, then $\alpha=\infty$ and by Equation \eqref{eq:detdieu}, we
have $\n(ps)=(\Det(\ga))^2=1$. Since $p,s\in\OOO$, we have
$\n(p)=\n(s)=1$ and $p,s\in\OOO^\times$, hence $\alpha,\beta$ are
Farey neighbours with respect to $\OOO$.  If $q\neq 0$, then
$p,q,r,s\in\OOO$, $\alpha= pq^{-1}$, $\beta= rs^{-1}$ and
$\n(qpq^{-1}s-qr)=(\Det(\ga))^2 =1$, hence $\alpha,\beta$ are Farey
neighbours with respect to $\OOO$.

Conversely, assume that $\alpha,\beta$ are Farey neighbours with
respect to $\OOO$. First assume that there exists $p,q,r,s\in \OOO$
such that $\alpha= pq^{-1}$, $\beta=rs^{-1}$, $q=0$ and $p,s\in
\OOO^\times$. Then $\ga=\begin{bsmallmatrix} p&r\\0&
s\end{bsmallmatrix}$ belongs to $\Ga_\OOO$, and $\alpha= \infty=
\ga\cdot \infty$ and $\beta=rs^{-1}=\ga\cdot 0$, as wanted.
Otherwise, there exists $p,q,r,s\in \OOO$ such that $\alpha= pq^{-1}$,
$\beta= rs^{-1}$, $q\neq0$ and $\n(qpq^{-1}s-qr)=1$. Then
$\begin{psmallmatrix} p&r\\q&s \end{psmallmatrix}$ belongs to
$\SL_2(\OOO)$ by the preliminary comment. Hence $\ga
=\begin{bsmallmatrix} p&r\\q&s\end{bsmallmatrix}$ belongs to
$\PSL_2(\OOO)$ and maps $\infty$ and $0$ to $\alpha$ and $\beta$
respectively, as wanted.
\cqfd

\medskip
\noindent{\bf Proof of Theorem \ref{theo:quaternionfareycount}.}  We
apply Theorem \ref{theo:quatercount} to $\Ga=\Ga_\OOO$ and $(x,y)=
(0,\infty)$, so that by Lemma \ref{lem:fromarithtogeomquat}, we have
$\fn_{\OOO}= \Ga\cdot\{x,y\}$. The locally geodesic line
$\ell_{0,\infty}$ in $\Ga_\OOO\bs\HH^5_\RR$ is reciprocal as in the
rational case (the order two element $\begin{bsmallmatrix} 0&-1\\1&\ 0
\end{bsmallmatrix} \in\Ga_\OOO$ exchanges the two points at infinity
of $\wt\ell_{0,\infty} (\RR)$), hence $\iota_{\Ga_\OOO,\rm rec}(0,\infty)
=1$. The pointwise stabiliser in $\Ga_\OOO$ of the geodesic line
$\wt\ell_{0,\infty}(\RR)$ has cardinality $\frac{|\OOO^\times|^2}2$,
hence $m_{\Ga_\OOO}(0,\infty)= \frac{|\OOO^\times|^2}2$. The index in
$\Ga_{\OOO,\infty}$ of its unipotent subgroup of translations by
$\OOO$ is equal to $\frac{|\OOO^\times|^2}2$. Hence replacing the
quotient modulo $\Ga_{\OOO,\infty}$ in the left-hand side of the
formula in Theorem \ref{theo:quatercount} by the quotient modulo
$\OOO$ amounts to multiplying the right-hand side by
$\frac{|\OOO^\times|^2}2$. Therefore Theorem
\ref{theo:quaternionfareycount} does follow from Theorem
\ref{theo:quatercount}. \cqfd

{\small \bibliography{../biblio} }

\begin{thebibliography}{BAPP}

\bibitem[AES2]{AkaEinSha16Inv}
M.~Aka, M.~Einsiedler and U.~Shapira.
\newblock {\it Integer points on spheres and their orthogonal lattices}.
\newblock {Invent. Math.  {\bf 206} (2016) 379--396}.

\bibitem[Art]{Artin91}
E.~Artin.
\newblock {\it Algebra}.
\newblock {2nd ed., Prentice Hall, 1991}.

\bibitem[Ath]{Athreya16}
J.~Athreya.
\newblock {\it Gap distributions and homogeneous dynamics}.
\newblock {In "Geometry, topology, and dynamics in negative curvature", 
pp 1--31, London Math. Soc. Lect. Note Ser. {\bf 425}, 
Cambridge Univ. Press, 2016}.

\bibitem[BeS]{BesSav12}
M.~Bestvina and G.~Savin.
\newblock {\it Geometry of integral binary Hermitian forms}.
\newblock {J. Algebra {\bf 360} (2012) 1--20}.

\bibitem[BS]{BocSis22}
F.~P.~Boca and M.~Siskaki.
\newblock {A note on the pair correlation of Farey fractions}.
\newblock {Acta Arith. {\bf 205} (2022) 121--135}.

\bibitem[BZ]{BocZah05}
  F.~P.~Boca and A.~Zaharescu.
\newblock {\it The correlations of Farey fractions}.
\newblock {J. Lond. Math. Soc. {\bf 72} (2005) 25--39}.

\bibitem[Coh]{Cohn80}
H.~Cohn.
\newblock {\it Advanced number theory}.
\newblock {Dover, 1980. Reprint of {\it A second course in number
      theory}, Wiley, 1962}.

\bibitem[Cox]{Cox13}
D.~A.~Cox.
\newblock {\it Primes of the form $x^2 + ny^2$}.
\newblock {2nd ed., Pure Appl. Math., Wiley, 2013}.

\bibitem[Cre]{Cremona84}
J.~E.~Cremona.
\newblock {\it Hyperbolic tessellations, modular symbols, and
  elliptic curves over complex quadratic fields}.
\newblock {Compos. Math. {\bf 51} (1984) 275--324}.

\bibitem[Duk]{Duke03}
W.~Duke.
\newblock {\it Rational Points on the Sphere}.
\newblock {Rankin memorial issues. Ramanujan J. {\bf 7} (2003) 235--239}.

\bibitem[EMV]{EllMicVen13}
J.~Ellenberg, P.~Michel and A.~Venkatesh.
\newblock {\it Linnik’s ergodic method and the distribution
of integer points on spheres}.
\newblock {Automorphic Rep. $L$-Functions {\bf 22} (2013),
119--185}.

\bibitem[EGM]{ElsGruMen98}
J.~Elstrodt, F.~Grunewald and J.~Mennicke.
\newblock {\it Groups acting on hyperbolic space: Harmonic
  analysis and number theory}.
\newblock {Springer Mono. Math., Springer Verlag, 1998}.

\bibitem[For]{Ford38}
L.~R.~Ford.
\newblock {\it Fractions}.
\newblock {Amer. Math. Monthly {\bf 45} (1938) 586--601}.

\bibitem[Hal]{Hall70}
R.~R.~Hall.
\newblock {\it A note on Farey series}.
\newblock {J. Lond. Math. Soc. {\bf 2} (1970) 139--148}.

\bibitem[HaT]{HalTen84}
R.~R.~Hall and G.~Tenenbaum.
\newblock {\it On consecutive Farey arcs}.
\newblock {Acta Arith. {\bf 44} (1984) 397–405}.

  
\bibitem[Hay]{Haynes04}
A.~Haynes.
\newblock {\it The distribution of special subsets of the Farey sequence}.
\newblock {J. Number Theory {\bf 107} (2004) 95--104}.

\bibitem[Hee]{Heersink16}
B. Heersink.
\newblock {\it Poincaré sections for the horocycle flow in
  covers of $\SL(2,\RR)/\SL_2(\ZZ)$ and applications to
  Farey fraction statistics}.
\newblock {Monatsh. Math. {\bf 179} (2016) 389--420}.

\bibitem[HK]{HorKar23}
  T.~Horesh and Y.~Karasik.
\newblock {\it Equidistribution of primitive lattices in $\RR^n$}.
\newblock {Q. J. Math. {\bf 74} (2023) 1253--129}.

\bibitem[HN]{HorNev23}
T.~Horesh and A.~Nevo.
\newblock {\it Horospherical coordinates of lattice points in 
hyperbolic space: effective counting and equidistribution}.
\newblock {Pacific J. Math. {\bf 324} (2023) 265–294}.

\bibitem[Kel]{Kellerhals03}
R.~Kellerhals.
\newblock {\it Quaternions and some global properties of 
hyperbolic $5$-manifolds}.
\newblock {Canad. J. Math. {\bf 55} (2003) 1080--1099}.

\bibitem[KO]{KraOse90}
V.~Krafft and D.~Osenberg.
\newblock {\it Eisensteinreihen f\"ur einige arithmetisch definierte
Untergruppen von $\SL_2(\HH)$}.
\newblock {Math. Z. {\bf 204} (1990) 425--449}.

\bibitem[Lem]{Lemmermeyer21}
F.~Lemmermeyer.
\newblock {\it Quadratic number fields}.
\newblock {Springer Verlag,  2021}.

\bibitem[Lut]{Lutsko22}
C.~Lutsko.
\newblock {\it Farey sequences for thin groups}.
\newblock {Int. Math. Res. Not.  {\bf 15} (2022) 11642--11689}.

\bibitem[Man1]{Manin72}
Y.~Manin.
\newblock {\it Parabolic points and zeta functions of modular curves}.
\newblock {Math. USSR Izvestija {\bf 6} (1972) 19--64}.

\bibitem[Man2]{Manin09}
Y.~Manin.
\newblock {\it Lectures on modular symbols}.
\newblock {In ``Arithmetic geometry'', pp. 137--152,
  Clay Math. Proc. {\bf 8}, Amer. Math. Soc. 2009}.

\bibitem[Mar1]{Marklof10Inv}
J.~Marklof.
\newblock {\it The asymptotic distribution of Frobenius numbers}.
\newblock {Invent. Math. {\bf  181}  (2010) 179--207}.

\bibitem[Mar2]{Marklof10AMS}
J.~Marklof.
\newblock {\it Horospheres and Farey fractions}.
\newblock {In "Dynamical numbers---interplay between dynamical 
systems and number theory",  Contemp. Math. {\bf 532} 97--106, 
Amer. Math. Soc. 2010}.

\bibitem[Mar3]{Marklof13}
J.~Marklof.
\newblock {\it Fine-scale statistics for the multidimensional Farey sequence}.
\newblock {In "Limit theorems in probability, statistics and number
theory", pp.~49–-57, P.~Eichelsbacher et al eds.
Springer Proc. Math. \& Stat. {\bf 42}, Springer Verlag, 2013}.

\bibitem[McM]{McMullen21}
C.~T.~McMullen.
\newblock {\it Modular symbols for Teichmüller curves}.
\newblock {J. Reine Angew. Math. {\bf 777} (2021) 89--125}.

\bibitem[Men]{Mendoza80}
E.~Mendoza.
\newblock {\it Cohomology of $\PGL_2$ over imaginary quadratic integers}.
\newblock {Bonner Math. Schriften {\bf 148}, Universität Bonn, 1980}.

\bibitem[PP1]{ParPau11BLMS}
J.~Parkkonen and F.~Paulin.
\newblock {\it On the representations of integers by 
indefinite binary Hermitian forms}.
\newblock {Bull. Lond. Math. Soc. {\bf 43} (2011) 1048--1058}.

\bibitem[PP2]{ParPau13ANT}
J.~Parkkonen and F.~Paulin.
\newblock {\it On the arithmetic and geometry of 
binary Hamiltonian forms}.
\newblock {Appendix by Vincent Emery. 
Algebra Number Theory {\bf 7} (2013) 75--115}.
  
\bibitem[PP3]{ParPau17ETDS}
J.~Parkkonen and F.~Paulin.
\newblock {\it Counting common perpendicular arcs 
in negative curvature}.
\newblock {Ergodic Theory Dynam. Systems {\bf 37} (2017) 900--938}.

\bibitem[PP4]{ParPau24}
J.~Parkkonen and F.~Paulin.
\newblock {\it Joint partial equidistribution of Farey rays in
negatively curved manifolds and trees}.
\newblock {Ergodic Theory Dynam. Systems {\bf 44} (2024) 2700-2736}.

\bibitem[PP5]{ParPau25Ingham}
J.~Parkkonen and F.~Paulin.
\newblock {\it Divergent geodesics, ambiguous closed geodesics and
  the binary additive divisor problem}.
\newblock {Preprint {\tt [arXiv:2409.18251]}}.

\bibitem[PPS]{ParPauSay25}
J.~Parkkonen, F.~Paulin and R.~Sayous.
\newblock {\it Equidistribution of divergent geodesics in negative
  curvature}.
\newblock {Preprint {\tt [arXiv:2501.03925]}}.

\bibitem[Sar]{Sarnak07}
P.~Sarnak.
\newblock {\it Reciprocal geodesics}.
\newblock {Clay Math. Proc. {\bf 7} (2007) 217--237}.

\bibitem[Say]{Sayous25}
R.~Sayous.
\newblock {\it Gaps in the complex Farey sequence of an
imaginary quadratic field}.
\newblock {Int. J. Number Theory {\bf 21} (2025) 2015--2036}.

\bibitem[Sch]{Schmidt68}
W.~Schmidt.
\newblock {\it Asymptotic formulae for point lattices of bounded
  determinant and subspaces of bounded height}.
\newblock {Duke Math. J. {\bf 35} (1968) 327--339}.

\bibitem[Shi]{Shimura71}
G.~Shimura.
\newblock {\it Introduction to the arithmetic theory of automorphic functions}.
\newblock {Princeton Univ. Press, 1971}.

\bibitem[Som]{Sommer11}
J. Sommer.
\newblock {\it Introduction à la théorie des nombres algebriques}.
\newblock {Hermann, 1911}.

\bibitem[Vig]{Vigneras80}
M.~F.~Vign{\'e}ras.
\newblock {\it Arithm\'etique des alg\`ebres de quaternions}.
\newblock {Lect. Notes in Math. {\bf 800}, Springer Verlag, 1980}.

\bibitem[Wu]{Wu02}
J.~Wu.
\newblock {\it On the Primitive Circle Problem}.
\newblock {Monatsh. Math. {\bf 135} (2002) 69--81}.

\bibitem[Zul]{Zullig28}
J.~Züllig.
\newblock {\it Geometrische Deutung unendlicher Kettenbrüche
und ihre Approximation durch rationale Zahlen}.
\newblock {Orell Füssli Verlag, 1928}.

\end{thebibliography}

\bigskip
{\small
\noindent \begin{tabular}{l} 
Department of Mathematics and Statistics, P.O. Box 35\\ 
40014 University of Jyv\"askyl\"a, FINLAND.\\
{\it e-mail: jouni.t.parkkonen@jyu.fi}
\end{tabular}

\medskip\noindent\begin{tabular}{l}
Laboratoire de mathématique d'Orsay, UMR 8628 CNRS,\\
Universit\'e Paris-Saclay, 91405 ORSAY Cedex, FRANCE\\
{\it e-mail: frederic.paulin@universite-paris-saclay.fr}
\end{tabular}
}

\end{document}